\documentclass[11pt,twoside]{amsart}
\usepackage{array,hhline,amssymb,amsthm,calrsfs,xspace}
\usepackage[dvips]{graphics}
\usepackage[vcentermath]{youngtab}
\usepackage{diagrams}
\diagramstyle{labelstyle=\scriptstyle,small}
\newarrow{Eqto}=====

\def\VEC#1,#2{(#1_1,#1_2,\dots,#1_{#2})}
\def\OVEC#1,#2{(#1_0,#1_1,#1_2,\dots,#1_{#2})}
\def\IVEC#1,#2,#3{({#1}_{{#2}_1},{#1}_{{#2}_2},\dots,{#1}_{{#2}_{#3}})}
\def\SET#1,#2{\{#1_1,#1_2,\dots,#1_{#2}\}}
\def\OSET#1,#2{\{#1_0,#1_1,\dots,#1_{#2}\}}
\def\ISET#1,#2,#3{\{{#1}_{{#2}_1},{#1}_{{#2}_2},\dots,{#1}_{{#2}_{#3}}\}}
\def\FAM#1,#2{#1_1,#1_2,\dots,#1_{#2}}
\def\OFAM#1,#2{\ #1_0,#1_1,\dots,#1_{#2}\ }
\def\IFAM#1,#2,#3{{#1}_{{#2}_1},{#1}_{{#2}_2},\dots,{#1}_{{#2}_{#3}}}
\def\WP#1,#2{#1_1\w #1_2\w\dots\w #1_{#2}}
\def\OWP#1,#2{#1_0\w #1_1\w\dots\w #1_{#2}}
\def\IWP#1,#2,#3{#1_{#2_1}\w #1_{#2_2}\w\dots\w #1_{{#2_{#3}}}}
\def\POL#1,#2,#3{#1_0#2^{#3}+#1_1#2^{#3-1}+\dots+#1_{#3-1}#2+#1_{#3}}
\def\LC#1,#2,#3{{#1}_1{#2}_1+{#1}_2{#2}_2+\dots+{#1}_{#3}#2_{#3}}

\renewcommand{\ge}{\geqslant}
\renewcommand{\le}{\leqslant}
\let\w=\wedge
\let\x=\times
\let\ox=\otimes
\let\PD=\partial
\let\n=\nabla

\newcommand{\EPbox}{\nopagebreak[4]\hfill$\square$\par}
\newcommand{\PDB}{\overline\PD}

\def\Ddots{\mathinner{\mkern1mu\raise1pt\hbox{.}\mkern2mu\raise5pt\hbox
                  {.}\mkern2mu\raise8pt\vbox{\kern7pt\hbox{.}}\mkern1mu}}
\def\bydef{\stackrel{\text{def}}{=}}
\def\ovl#1{\overline{#1}}
\def\wtd#1{\widetilde{#1}}
\def\wht#1{\widehat{#1}}
\def\({\left(}
\def\){\right)}
\def\OP{\mathop\oplus\limits}
\def\X{\mathop\x\limits}
\def\OX{\mathop\ox\limits}
\def\comp{\mathop{\scriptstyle\circ}}
\def\bdot{{\scriptscriptstyle\bullet}}

\def\bl{\!\in\!}
\def\fain#1#2{\;\forall\,{#1}\bl{#2}\;}
\def\<#1,#2>{\left\langle\,#1\,,\,#2\,\right\rangle}
\newcommand{\rep}[1]{\left\langle{#1}\right\rangle}
\let\IFF=\Longleftrightarrow
\let\THEN=\Longrightarrow

\def\Ann{\mathrm{Ann}\,}

\def\codim{\mathrm{codim}\,}

\def\End{\mathrm{End}}

\def\ext{\mathrm{Ext}}

\def\GL{\mathrm{GL}}
\def\GLc{$\mathrm{GL}$-cha\-rac\-ter\xspace}

\def\Gr{\mathrm{ Gr}}

\def\hom{\mathrm{Hom}}
\def\Id{\mathrm{Id}}

\def\im{\mathrm{im}\,}
\def\Lie{\cL\!\mathit{ie}}

\def\Pic{\mathrm{Pic}}

\def\Skew{\mathrm{Skew}}
\def\SL{\mathrm{SL}}
\def\SO{\mathrm{SO}}

\def\Stab{\mathrm{Stab}}

\def\SU{\mathrm{SU}}
\def\Sym{\mathrm{Sym}}

\def\Tor{\mathrm{Tor}}
\def\tr{\mathrm{tr}\,}
\def\Tr{\mathrm{Tr}\,}

\def\ga{\mathfrak{a}}
\def\gb{\mathfrak{b}}

\let\GG=\gg
\def\gg{\mathfrak{g}}
\def\gh{\mathfrak{h}}

\def\gp{\mathfrak{p}}

\def\gC{\mathfrak{C}}

\def\gS{\mathfrak{S}}

\newcommand{\fA}{\ensuremath{\mathsf A}}
\newcommand{\fB}{\ensuremath{\mathsf B}}
\newcommand{\fL}{\ensuremath{\mathsf L}}

\newcommand{\cA}{\ensuremath{\mathcal A}}

\newcommand{\cE}{\ensuremath{\mathcal E}}

\newcommand{\cH}{\ensuremath{\mathcal H}}

\newcommand{\cL}{\ensuremath{\mathcal L}}

\newcommand{\cO}{\ensuremath{\mathcal O}}

\def\G{\Gamma}
\def\L{\Lambda}
\let\O=\Omega

\let\a=\alpha
\let\b=\beta
\let\f=\varphi
\let\g=\gamma
\let\d=\delta
\let\D=\Delta
\let\e=\varepsilon
\let\k=\varkappa
\let\l=\lambda
\let\o=\omega
\let\r=\varrho
\let\s=\sigma
\let\t=\vartheta

\def\CC{\mathbb C}

\def\NN{\mathbb N}
\def\PP{\mathbb P}

\def\TT{\mathbb T}
\def\ZZ{\mathbb Z}

\newcommand{\Lt}{L_{\ge2}}
\newcommand{\without}[2]
    {\mathop{\wht{\hphantom{{#1}_{#2}}
                  \vphantom{{#1}_{#2}}}}\limits_{#2}}

\newcommand{\wrt}{w.\,r.\,t. }
\newcommand{\etc}{etc.}
\newcommand{\ie}{i.\,e. }

\newcommand{\SHW}{SHW\--or\-bit\xspace}
\newcommand{\SHWs}{SHW\--or\-bits\xspace}
\newcommand{\hw}{\mathrm{hw}}
\newcommand{\hqt}{h-\-tab\-leau}
\newcommand{\Dol}{Dolbeault\xspace}
\newcommand{\LRr}{the Littlewood--Richardson rule\xspace}
\newcommand{\wl}{\L^{\mathrm{wt}}}
\newcommand{\rl}{\L^{\mathrm{rt}}}
\newcommand{\Lb}{\ovl L}
\newcommand{\dk}{d_{\text{\sc k}}}
\newcommand{\dc}{d_{\text{\sc c}}}

\newtheorem{theorem}[subsubsection]{Theorem}
\newtheorem{proposition}[subsubsection]{Proposition}
\newtheorem{lemma}[subsubsection]{Lemma}
\newtheorem{corollary}[subsubsection]{Corollary}

\theoremstyle{definition}
\newtheorem{example}[subsubsection]{Example}

\theoremstyle{remark}

\numberwithin{equation}{section}

\date{July 28, 2006}

\begin{document}



\title{On syzygies of highest weight orbits}
\thanks{the first author thanks scientific school supporting grant NSh-8004.2006.2 and the
second author~--- RFBR grant 04-02-16538 for partial support}

\author{A.~L.~Gorodentsev}
\address{Institute of Theoretical and Experimental Physics (Moscow)}
\curraddr{B. Cheremushkinskaya, 25, 117259 Moscow, Russia}
\email{gorod@itep.ru}

\author{A.~S.~Khoroshkin}
\address{Institute of Theoretical and Experimental Physics (Moscow)}
\curraddr{B. Cheremushkinskaya, 25, 117259 Moscow, Russia}
\email{khorosh@itep.ru}

\author{A.~N.~Rudakov}
\address{Norwegian University of Science and Technology (Trondheim)}
\curraddr{Math. Dept. NTNU,
          Alfred Getz vei 1, NO-7034 Trondheim, Norway}
\email{rudakov@math.ntnu.no}

\begin{abstract}
  We consider the graded space $R$ of syzygies for the coordinate algebra $A$ of projective variety $X=G/P$ embedded into projective space as an orbit of the highest weight vector of an irreducible representation of semisimple complex Lie group $G$. We show that $R$ is isomorphic to the Lie algebra cohomology $H=H^\bdot(\Lt,\CC)$, where $\Lt$ is graded Lie subalgebra of the graded Lie s-algebra $L=L_1\oplus\Lt$ Koszul dual to $A$. We prove that the isomorphism identifies the natural associative algebra structures on $R$ and $H$ coming from their Koszul and Chevalley DGA resolutions respectively. For subcanonically embedded $X$ a Frobenius algebra structure on the syzygies is constructed. We illustrate the results by several examples including the computation of syzygies for the Pl\"ucker embeddings of grassmannians $\Gr(2,N)$.
\end{abstract}

\maketitle

\section*{Introduction}

   A.~Losev brought to our attention the fact that some computations made by N.~Ber\-ko\-vits  in the framework of string theory\footnote{see \cite{Bk1}, \cite{Bk2}} contain an intricate description of minimal resolution for the projective coordinate algebra of connected component of complex isotropic grassmannian $\Gr^+_{\mathrm{iso}}(5,10)$ associated with non degenerate quadratic form on $\CC^{10}$. Subsequent papers by M.~Movshev, A.~Schwarz and others\footnote{see \cite{CN}, \cite{Mo1}, \cite{Mo2}, \cite{MS1}, \cite{MS2}, \etc} clarify the Ber\-ko\-vits computations as well as the interplay between the syzygies and a graded Lie superalgebra Koszul dual\footnote{in the sense of \cite{GK}} to the coordinate algebra of $\Gr^+_{\mathrm{iso}}(5,10)$. But the main focus of these papers remains with s-symmetric field theories, and the presentation of underlying mathematics seems to us rather tangled and over-diligent.

   The main goal of these notes is to give relatively simple and clear presentation of the mathematics behind the sophisticated computations cited above and to illustrate it by some classical geometric examples. We tried to make the text self contained and understandable for advanced students. We believe that the subject is fruitful, and hope that this exposition could promote mutual understanding between mathematical and physical communities.

   The paper consists of five sections. The first three of them are essentially independent.  In \S1 we collect necessary algebraic geometric properties of projective varieties
   $$X=\PP(G\cdot v_{\mathrm hw})\subset\PP(V)
   $$
   that appear as the projectivization of the orbit of the highest
   weight vector in an irreducible representation $V$ of a simply
   connected semisimple complex Lie group $G$. Class of these
   varieties clearly contains the isotropic grassmannians and
   serves a convenient generalization for the smooth\footnote{In
   fact, besides the grassmannians, Berkovits and others
   consider also a singular quadratic cone over
   $\Gr_{\mathrm{iso}}(5,11)$, which is especially interesting for
   physicists because it serves some model of gravity. But the
   mathematical statements here are rather unclear to us and most
   likely require appropriate $A_\infty$ enchantment of the smooth
   case. In these notes we restrict ourselves to the smooth
   homogeneous spaces only.} varieties of `pure spinors' considered
   in s-symmetric theories cited above. In \ref{Qrel} we write down explicitly the quadratic equations generating the homogeneous ideal\footnote{these equations go back to Kostant and are known to specialists (comp. with \cite{LT}, \cite{Li}); they have famous
   infinite dimensional extension developed by
   V.~Kac and D.~Peterson (see \cite{KP1}, \cite{KP2}); for convenience of
   readers we sketch a short geometric proof for the finite dimensional case}
   $$I(X)=\{f\in S(V^*)\,|\;f|_X\equiv0\}\;.
   $$
   In \ref{HVprop} we show that if $X$ is subcanonical, \ie $\o_X=\cO_X(-N)$ for some $N\in\NN$, then non zero cohomologies $H^q(X,\cO_X(m))$ can appear only for $q=0$ or for $q=\dim X$.

   In \S2 we consider an arbitrary smooth\footnote{in fact all constructions and results of \S2 hold for locally complete intersection varies as well} subcanonical projective variety satisfying the previous vanishing conditions on $H^q(X,\cO_X(m))$. By \S1, all subcanonical projective highest weight orbits belong to this class. Using Movshev's strategy from \cite{Mo1}, \cite{MS2}, we show that for such $X$ the space $R$ of the syzygies for the graded projective coordinate algebra
   $$A=S(V^*)/I(X)=\OP_{m\ge0}H^0\(X,\cO_X(m)\)
   $$
   inherits a natural structure of the Frobenius algebra. In other
   words, we construct an s-commutative multiplication on $A$ and a
   trace form
   $$A\rTo^\tr\CC
   $$
   such that the bilinear pairing $(a,b)=\tr(a\cdot b)$ is non
   degenerated. This generalizes and clarifies the duality
   isomorphisms constructed in \cite{Mo1}, \cite{MS2} for the
   smooth varieties of pure spinors.

   In \S3 we deal with an arbitrary commutative Koszul quadratic
   algebra $A$. All the coordinate algebras of projective highest
   weight orbits are of this sort due to Bezrukavnikov's result
   \cite{Bz}. It is well known that such an algebra is canonically
   identified with the cohomology algebra
   $$A\simeq H^\bdot(L,\CC)
   $$
   of the graded Lie superalgebra $L=\OP_{m\ge1}L_m$ whose universal enveloping algebra $A^!=U(L)$ is Priddy dual to $A$. We show that there is an isomorphism of algebras
   $$R\simeq H^\bdot(\Lt,\CC)\;,
   $$
   where the space of syzygies $R$ is considered with the algebra
   structure constructed in \S2, and $\Lt=\OP_{m\ge2}L_m\subset L$ is
   graded Lie subalgebra started with $L_2$-component of $L$ (the
   algebra structure on $H^\bdot(\Lt)$ is standard). The proof is based
   on a kind of differential perturbation lemma (see \ref{dpl}).

   In \S4 we illustrate the previous technique by several non trivial
   examples. Namely, we compute the syzygies of the most singular
   commutative quadratic algebra $A=\TT(V^*)/\Skew(V^*\ox V^*)$ (see
   \ref{syz_free_Lie}), the syzygies of rational normal curves (see
   \ref{cor:VCsyz}), and the syzygies of grassmannians $\Gr(2,N)$ under
   the Pl\"ucker embeddings (see \ref{syzygies}). These results are
   also known for experts and can be extracted from \cite{Gr},
   \cite{JPW}, \cite{PW}, \cite{W} and references therein. Our approach
   allows to treat all three examples uniformly: we use the description
   of the Lie algebra cohomology $H^\bdot(\Lt)$. In the first two
   examples (actually served by free Lie algebras) the computation is
   very simple and takes just a few rows. In the grassmannian case the
   algebra of syzygies is what we call {\it a hook algebra\/}. In the
   last \S5 we collect generic properties of the hook algebras, in
   particular, we prove that each hook algebra is quadratic and koszul.

   \smallskip

   \subsection*{Acknowledgments.}
   We thank the Mittag-Lefler Institute, the Max Planck Institute, and
   IHES for the excellent possibilities to meet our colleagues and work
   on the subject, and NTNU for the nice opportunity to finish this
   job. We are grateful to B.~Feigin, V.~Kac, A.~Losev, M.~Movshev and
   D.~Piontkovsky for many useful discussions. We thank D.~Panyushev
   and V.~Popov for the references on the papers \cite{LT}, \cite{Li},
   G.Olshansky for the reference on \cite{IW}, and J.~Weyman for
   the references on \cite{JPW}, \cite{PW}, and \cite{W}.

\section{Projective orbit of the highest weight vector}

\subsection{Basic notations}\label{QHWcond}
  Let $G$ be connected and simply connected complex semisimple
  algebraic group and $V=V_\l$ be its complex irreducible linear
  representation with a highest weight $\l$. We fix Cartan and
  Borel subgroups $T\subset B\subset G$ and write $v_\hw\in V$ for
  a highest weight vector and $P\subset G$ for a parabolic subgroup
  stabilizing 1-dimensional subspace $\CC\cdot v_\hw$. In this
  part we consider the projectivization of the highest weight
  orbit
  \begin{equation}\label{cpe}
      X=G/P\simeq\PP(G\cdot v_\hw)\rInto^{\quad\f\quad}\PP(V)\;,
  \end{equation}
  which is a homogeneous $G$-space with the natural left action of
  $G$. We always put $\dim X=d$, $\dim V=n+1$. Our especial
  interest is in the case when the canonical class of $X$ is a
  negative integer multiple of the hyperplane section, \ie
  \begin{equation}\label{QHWeq}
     \o_X=\cO_X(-N)\quad\text{for some }N\in\NN\;,
  \end{equation}
  where $\cO_X(1)=\f^*\cO_{\PP(V)}(1)$ is the ample line bundle
  coming from the projective embedding \eqref{cpe}. We will call
  such a variety $X$ {\it a subcanonical\/} highest weight orbit (or
  {\it a \SHW\/} for shortness). To clarify this condition, let us
  fix some standard notations related to Lie algebras and recall
  some basic facts about vector bundles on homogeneous spaces.

\subsubsection{Lie algebra notations.}
 We denote by $\gg\supset\gp\supset\gb\supset\gh$ the Lie algebras
 of $G\supset P\supset B\supset T$ and write
 $\rl\subset\wl\subset\gh^*$ for the root and weight lattices of
 $\gg$ and $\D\subset\rl$ for the set of all {\it positive\/} roots. As
 usual, we put
 $$\r=\frac12\,\sum_{\a\in\D}\a
 $$
 to be the half sum of all positive roots. Let $\{\a_i\}$ be the
 basic simple positive roots\footnote{recall that $(\r,\a_i)=1$ for
 any simple root $\a_i$} and $\{\o_i\}$ be the corresponding
 fundamental weights, which are dual to $\a_i$ \wrt the Killing
 form. Then we have
  \begin{align}
     \label{la-dec}
     \l&=\sum n_i\cdot\o_i\quad\text{with integer }n_i\ge0\,,\\
     \label{lga-dec}
     \gp&=\gb\oplus\bigoplus_{\a\in\D^\gp}\gg_{-\a}\,,\quad\text{where }
     \D_\gp=\D\cap\l^\bot
  \end{align}
  (a simple root $\a_i\in\D_\gp$ iff the corresponding $n_i=0$). We write
  $$\wl_\gp=\D_\gp^\bot\cap\wl=
    \{\,\mu\in\wl\,|\;(\mu,\a)=0\;\fain\a{\D_\gp}\}=
    \!\bigoplus_{i\,|\,n_i\ne0}\ZZ\cdot\o_i
  $$
  for the set of all weights that produce the characters for $\gp$ and
  write $\rep\mu$ for the 1-dimensional $P$-module coming from such
  a character $\mu\in\wl_\gp$. We also put
  $$\r_\gp=\frac12\,\sum_{\a\in\D_\gp}\a\;.
  $$
  Besides $B$, we will sometimes consider its opposite Borel
  subgroup $B'$. If $B=T\ltimes U$, where $U$ is the
  unipotent part of $B$, then $B'=T\ltimes U'$ and $U'TU$ is a
  dense open subset in $G$. Similarly, we will write $\r'=-\r$ for
  the half sum of all the negative roots.

\subsubsection{Vector bundles on $G/P$.}\label{vb-sec}
  With any representation $E$ of $P$ is associated a vector
  bundle
  $$\wtd E=G\X_PE
  $$
  over $X=G/P$ with the fiber $E$. Its total space consists of
  pairs $(g,e)\in G\x E$ modulo the equivalence $(gp,e)\sim(g,pe)$
  for $p\in P$. Global sections $X\rInto^{\;s_f}\wtd E$ are
  naturally identified with functions $G\rTo^f E$ such that
  \begin{equation}\label{sec-fun}
        pf(g)=f(gp^{-1})\quad\text{for all }p\in P\;.
  \end{equation}
  The left $G$-action on $G/P$ is extended canonically to the left
  $G$-action on the whole of $\wtd E$ by the rule
  $$g\cdot(g_1,e)=(gg_1,e)
  $$
  and induces a linear representation of $G$ in the space
  $\G(X,\wtd E)$, of the global sections of $\wtd E$. In terms of
  equivariant functions \eqref{sec-fun}, in this representation
  an element $g\in G$ sends a function $f$ to a function $g\cdot f$
  defined by prescription
  \begin{equation}\label{fun-act}
     g\cdot f(g_1)=f(g^{-1}g_1)\;.
  \end{equation}
  If the action of $P$ on $E$ can be extended to an
  action of $G$, then we have a vector bundle isomorphism
  $\wtd E\rTo^\sim X\x E$, which takes an equivalence class of
  $(g,e)$ to $(g,ge)$.

  In particular, each 1-dimensional $P$-module $\rep\mu$\,, where
  $\mu\in\wl_\gp$\,, leads to the line bundle
  \begin{equation}\label{albundef}
     [\mu]=G\X_P\rep\mu\;.
  \end{equation}
  We will write $[\mu]$ for this bundle considered as an element of
  the Picard group\footnote{in fact, $\Pic(X)$ is spanned by the
  line bundles \eqref{albundef}, see \cite{VP}} $\Pic(X)$ and will write
  $\cO_X\rep\mu$ for the corresponding invertible sheaf of its
  local sections.

  For example, 1-dimensional $P$-module $\rep\l=\CC\cdot
  v_\hw\subset V_\l$, spanned by the highest weight vector in
  $G$-module $V_\l$, produces the tautological line
  subbundle
  $$\cO_X\rep\l=\cO_X(-1)\subset\wtd V_\l=X\x V_\l\;,
  $$
  which coincides with the restriction of the
  tautological line bundle $\cO_{\PP(V_\l)}(-1)$ onto $X$.

   It is easy to see that the tangent bundle $T_X=G\X_PQ$ comes from
   the representation
   $$Q=\gg / \gp=\!\!\!\!\!
      \bigoplus_{\a\in\D\smallsetminus\D_\gp}\!\!\!\!\gg_{-\a}\;.
   $$
  Hence, the anticanonical line bundle is expressed in $\Pic(X)$
  as the sum
  $$\o^*_X=\L^dT_X=\!\!\!\!\!\!\sum_{\a\in\D\smallsetminus\D_\gp}[-\a]=
          -2\,(\r-\r_\gp)\;.
  $$
  Thus, a projective orbit \eqref{cpe} is subcanonical in
  $\PP(V_\l)$, \ie satisfies \eqref{QHWeq}, iff
  \begin{equation}\label{SHWcond-weights}
      2\,(\r-\r_\gp)=N\,\l\quad\text{for some }N\in\NN\;.
  \end{equation}
  This is significant restriction on $\l$ (e.\,g. we will see in
  \ref{H0ofOn} that it forces quite strong vanishing condition on
  the cohomologies of invertible sheaves $\cO_X(k)$).

  In the next
  examples we use the standard Bourbaki notations from
  \cite{Bu2}.

\begin{example}[the grassmannian $\Gr(2,5)$]\label{g25}
  Let $G=\SL(5,\CC)$ with the diagonal torus $T\subset G$,
  $$\gh=\{\LC a,\e,5\,|\;{\textstyle \sum}\,a_i=0\}\;,
  $$
  and the simple roots\;
  $\a_i=\e_i-\e_{i+1}$\,, $1\le i\le4$. Then the projective embedding
  $$\Gr(2,5)\rTo^{\;\f_\l\;}_\sim X=\PP(G\cdot v_\hw)\subset\PP(V_\l)
  $$
  by means of the representation $V_\l$ with the highest weight
  $\l=[0,m,0,0]=m\cdot\o_2$\,, $m\ge1$\,,
  is subcanonical only for $m=1,\,5$. Indeed,  we have
  \begin{gather*}
          \o_2=\(3\,\e_1+3\,\e_2-2\,\e_3-2\,\e_4-2\,\e_5\)/5\;,\\
         2\,\r=4\,\e_1+2\,\e_2-2\,\e_4-4\,\e_5\;,\\
     2\,\r_\gp=\e_1-\e_2+2\,\e_3-2\,\e_5\;.
  \end{gather*}
  Thus, $\l=m\o_2$ divides $2\,(\r-\r_\gp)=5\,\o_2$ only for $m=1$
  and $m=5$. The former, non tautological, case gives the Pl\"ucker
  embedding $\Gr(2,5)\rInto\PP_9$ and leads to $N=5$ in
  \eqref{QHWeq}.
\end{example}

\begin{example}[the group $G=\SL(3,\CC)\x\SL(2,\CC)$]\label{SL3xSL2}
  This group\footnote{note that it coincides with the semisimple component for the complexification of the
  compact Lie group $\SU(3)\x\SU(2)\x\SU(1)$}
  has the diagonal torus $T\subset G$ with $\gh=\gh'\oplus\gh''$, where
  \begin{gather*}
        \gh'=\{a_1\e_1'+a_2\e_2'+a_3\e_3'\,|\;
                a_1+a_2+a_3=0\}\;,\\
       \gh''=\{b_1\e_1''+b_2\e_2''\,|\;b_1+b_2=0\}\;,\\
       2\,\r=2\,\e_1'-2\,\e_3'+\e_1''-\e_2''\;.
    \end{gather*}
  The representation $V_\l=V'\ox V''$ (where $V'$, $V''$ are the tautological 3- and 2-dimensional $\SL(3)$- and $\SL(2)$-modules) has
  $$\l=\o_1'+\o_1''=
    \(2\,\e_1'-\e_2'-\e_3'\)/3+\(\e_1''-\e_2''\)/2\;.
  $$
  For $P=\Stab(v_\l)$ we have $2\r_\gp=\e_2'-\e_3'$, which gives
  $$2\,(\r-\r_\gp)=2\,\e_1'-\e_2'-\e_3'+\e_1''-\e_2''=3\,\o_1'+2\,\o_2''\;.
  $$
  We conclude that the tautological highest weight embedding of
  $G/P$ is itself not subcanonical. But the embedding
  $$G/P\rTo^{\;\f_\mu\;}_\sim X=\PP(G\cdot v_\mu)\subset
    \PP(S^3V'\ox S^2V'')
  $$
  (corresponding to $V_\mu=S^3V'\ox S^2V''$ with $\mu=3\o_1'+2\o_2''$)
  is subcanonical and \eqref{QHWeq} holds with $N=1$.
\end{example}

\begin{example}[even dimensional pure spinors]\label{so10}
  Let $G=\mathrm{Spin}(10,\CC)$ be the universal covering for
  $\SO(10,\CC)$ and $Y=\Gr^+_{\mathrm{iso}}(5,10)$ be a connected
  component of the grassmannian of 5-dimensional isotropic
  subspaces\footnote{since $Y$ does also parameterize the
  decomposable elements of the Clifford algebra (see \cite{Ch}), it
  is often called the variety {\it of 10-dimensional pure spinors\/}} in
  $\CC^{10}$ (this is the case originally considered by Berkovits
  in \cite{Bk1}). Here $\gg$ is the semisimple Lie algebra of type
  $D_5$. Using the standard notations of Bourbaki (see \cite{Bu2})
  as above, we can write
  $$2\,\r=8\,\e_1+6\,\e_2+4\,\e_3+2\,\e_4
  $$
  and $Y=G/P$, where $P$ has
  $2\,\r_\gp=4\,\e_1+2\,\e_2-2\,\e_4-4\,\e_5$\,.  This gives
  $$2\,(\r-\r_\gp)=4\,(\e_1+\e_2+\e_3+\e_4+\e_5)\;.
  $$
  The HW-orbit embedding corresponding to
  $\l=\o_5=\(\e_1+\e_2+\dots+\e_5\)/2$
  $$Y\rTo^{\;\f_{\o_5}\;}_\sim\PP(G\cdot v_{\o_5})
    \subset\PP(V_{\o_5})\;,
  $$
  is subcanonical with $N=8$. More generally, for any $m$ the
  variety of $2m$-di\-men\-si\-on\-al pure spinors
  $\Gr^+_{\mathrm{iso}}(m,2m)$ has subcanonical
  HW-embedding into $\PP(V_\l)$ with
  $\l=\(\e_1+\e_2+\dots+\e_m\)/2$\,.
  Indeed,
  \begin{gather*}
      2\,\r=2\,(m-1)\,\e_1+2\,(m-2)\,\e_2+\dots+2\,\e_{m-1}\\
  2\,\r_\gp=(m-1)\,\e_1+(m-3)\,\e_2+\dots-(m-1)\,\e_m\;,
    \end{gather*}
  and \eqref{QHWeq} holds with $N=2\,(m-1)$.
\end{example}

\subsection{Cohomologies of line bundles}\label{H0ofOn}
  The computation of cohomologies of line bundles on
  $X=G/P$ is reduced to the computation of cohomologies on the flag
  variety $Y=G/B$ via $G$-equivariant projection $Y\rTo^\pi X$\,.
  Namely, for each weight $\mu\in\wl_\gp$ we can consider the
  restriction $\rep\mu_B$, of the 1-dimensional $P$-module
  $\rep\mu$ onto $B$, and form a line bundle
  $[\mu]_B=G\X_B\rep\mu_B$, which is clearly isomorphic to the pull
  back of $[\mu]$ along $\pi$, \ie
  $\pi^*\cO_X\rep\mu=\cO_Y\rep\mu_B$. Then, the Leray spectral
  sequence gives canonical isomorphisms
  $$H^q(X,\cO_X\rep\mu)\simeq H^q(Y,\cO_Y\rep\mu_B)
  $$
  for all $q$. By this reason in the rest of this section we
  replace $P$ by $B$, $X$ by $Y$ and write simply $[\mu]$ and
  $\rep\mu$ instead of $[\mu]_B$, $\rep\mu_B$.

  It is convenient to describe the representation of $G$
  in the space $\G(Y,\cO_Y\rep\mu)$ in terms of its {\it lowest\/}
  vector. Namely, the lowest weight of $T$ on $\G(Y,\cO_Y\rep\mu)$
  is $\mu$ and a lowest weight vector is unique up to
  proportionality. Indeed, if we interpret the sections of
  $\cO_Y\rep\mu$ as the functions
  $$G\rTo^f\rep\mu\simeq\CC
  $$
  via \eqref{sec-fun} and \eqref{fun-act}, then for a lowest weight
  function $f$ we have $u'\cdot f=f$ for all $u'\in U'$ (see
  notations on page \pageref{vb-sec}). Hence, for any $t\in T$, $b\in
  B$
  $$t\cdot f(u'b)=f(t^{-1}u'b)=f(u''t^{-1}b)=
     f(t^{-1}b)\;,
  $$
  where $u''=t^{-1}ut\in U'$ fixes $f$ and $t^{-1}b\in B$. On the
  other hand, all three pairs
  $$(\,t^{-1}b\,,\,f(t^{-1}b)\,)\sim
    (\,e\,,\,\mu(t^{-1}b)(f(t^{-1}b))\,)\sim(\,e\,,\,f(e)\,)
  $$
  represent the same point in the total space of the bundle
  $\cO_Y\rep\mu=G\X_B\rep\mu$. Hence,
  $$f(t^{-1}b)=\mu(t^{-1}b)^{-1}f(e)=\mu(t)f(b)
  $$
  (were we write $\mu(b)$ for the operator corresponding to $b\in
  B$ in the representation $\rep\mu$). This means that $t\cdot
  f=\mu(t)f$ over the open dense subset $U'B\subset G$ and,
  moreover, $f$ is uniquely defined there by its value $f(e)$. So,
  the weight of $f$ is $\mu$ and this weight subspace is
  1-dimensional.

\subsubsection{The Borel--Weil--Bott theorem,}\label{BWBth}
  being formulated in our notations, describes all
  $G$-modules $H^q(Y,\cO_Y\rep\mu)$ as irreducible
  representations presented by their {\it lowest\/}
  weights\footnote{the formulation most commonly used in the
  representation theory (see, for example, \cite{Ja}) actually
  describes $G$-modules $H^q(\cO\rep\mu)$ in terms highest weights
  but the underlying homogeneous space is always taken to be
  $G/B'$, that is the lowest weight vector orbit}. Namely, given
  $\mu\in\wl$, we consider its shift $\mu+\r'$, by the half sum all
  the {\it negative\/} roots. There are two possibilities:
  \begin{enumerate}
     \item
     $\mu+\r'$ lies in the interior part of some Weyl chamber $C$
     \item
     $\mu+\r'$ belongs to a wall separating the Weyl chambers
  \end{enumerate}
  In the first case $(\a,\mu+\r')\ne0\;\fain\a\D$ and there exist a
  unique weight $\mu'$ in the {\it lowest\/} Weyl chamber
  $C_{\mathrm{low}}$ and a unique element $w$ of the Weyl group
  such that
  \begin{equation}\label{BWBeq}
     \mu+\r'=w(\mu'+\r')
  \end{equation}
  (indeed, $w$ has to be the symmetry that takes $C_{\mathrm{low}}$
  to $C$ and then $\mu'$ is determined uniquely). In this case
  $H^q(Y,\cO_Y\rep\mu)\ne0$ iff $q$ equals the length\footnote{here
  the length should be defined \wrt the reflections by the walls of
  the lowest chamber $C_{\mathrm{low}}$} of $w$ in the Weyl
  group. This non zero space is an irreducible $G$-module of the
  {\it lowest\/} weight $\mu'$.

  In the second case $\mu+\r'$ is orthogonal to some root
  $\a\in\D$ and the equation \eqref{BWBeq} is unsolvable in a sense
  that $\mu+\r'$ is not congruent to any weight in the interior part
  of $C_{\mathrm{low}}$ modulo the Weyl group action. In this case
  $H^q(Y,\cO_Y\rep\mu)=0$ for all $q$.

  For example, the above description of $\G(Y,\cO_Y\rep\mu)$ fits
  the Borel--Weil--Bott setup as the case when $\mu'=\mu$,
  $w=e$, $q=0$. In particular, $\G(Y,\cO_Y\rep\mu)\ne0$ iff $\mu$
  lies in the lowest weights Weyl chamber $C_{\mathrm{low}}$.

\begin{proposition}\label{HVprop}
  If $X=G/P$ is a $d$-dimensional \SHW with $\o_X=\cO_X(-N)$, then
  all the non zero cohomologies $H^q\(X,\cO_X(k)\)$ are only
  $$H^0\(X,\cO_X(m)\)\quad\text{ and }\quad
    H^d\(X,\cO_X(-N-m)\)\;,
  $$ where $m\ge0$ in the both cases.
\end{proposition}

\begin{proof}
  Since $\l$ is a highest weight, $-m\l\in C_{\mathrm{low}}$ lies in
  the lowest weights chamber for all $m\ge0$. So, for all
  $\cO_X(m)=\cO_X\rep{-m\l}$ we have $H^0\(X,\cO_X(m)\)\ne0$ and
  $H^q\(X,\cO_X(m)\)=0$ when $q>0$. By the Serre duality, this
  implies
  $$H^d\(X,\cO_X(-N-m)\)=H^0\(X,\cO_X(m)\)^*\ne0
  $$
  and vanishing of all $H^q\(X,\cO_X(-N-m)\)$ with $q\ne d$.

  To manage the remaining values $m=1,\,2,\,\dots\,,\,(m-1)$, let
  us note that in the Borel--Weil--Bott setup (see \ref{BWBth})
  the triviality of the representation
  $$H^d\(X,\cO_X(-N)\)=H^0\(X,\cO_X\)^*=\CC
  $$
  means that $N\l+\r'=w(\r')$ for some $w$ from the Weyl group. So, all
  the weights
  \begin{equation}\label{intweights}
       \l+\r'\,,\;2\l+\r'\,,\;\ldots\;,\;(N-1)\l+\r'
  \end{equation}
  are internal points of the segment $I=\{\,x\l+\r'\,|\;0\le x\le
  N\,\}$ whose endpoints are $\r'$ and $N\l+\r'=w(\r')$. The both
  endpoints have Euclidean length $||\r'||=\sqrt{(\r',\r')}$, which
  is the minimal length of the weights lying in the interior part
  of a Weyl chamber. By the convexity arguments, all the interior
  points of $I$ lay strictly closer to the origin. So, the weights
  \eqref{intweights} can not be interior points of a chamber. Hence
  they are not congruent to the interior points of
  $C_{\mathrm{low}}$ modulo the Weyl group action and we deal with
  the second case of the Borel--Weil--Bott theorem (in the sense of
  \ref{BWBth}). Therefore, all the cohomologies
  $H^q\(X,\cO\rep{-m\l}\)$ do vanish for $1\le m\le(N-1)$.
\end{proof}

\subsection{\mathversion{bold}Quadratic equations for $X$}\label{Qrel}
  In this section we write explicit quadratic equations generating the homogeneous ideal of the projectivization of an arbitrary highest weight vector orbit \eqref{cpe} (not necessary subcanonical).  Infinite dimensional versions of two propositions below were proved by Kac and Peterson in \cite{KP1}, \cite{KP2}. Finite dimensional case goes back to Kostant (comp. with \cite{LT}, \cite{Li}). For the convenience of readers we sketch here an easy finite dimensional proof.

  Let $U(\gg)$ be the universal enveloping algebra of $\gg$.
  Consider the Casimir element
  \begin{equation}\label{Cas}
     \O=\sum a_ib_i\in U(\gg)\;,
  \end{equation}
  where $a_i$ and $b_i$ form a pair of dual bases of $\gg$ \wrt
  the Killing form. It does not depend on the choice of dual bases,
  lays in the center of $U(\gg)$, and acts on each irreducible
  $\gg$-module by a scalar operator:
  $$\O|_{V_\l}=c_\l\Id_{V_\l}\;.
  $$
  The constant $c_\l$ can be computed by the following formula
  (see, for example, \cite{FH})
  \begin{equation}\label{CasVal}
     c_\l=(\,\l+\r\,,\,\l+\r\,)-(\,\r\,,\,\r\,)
         =(\,\l+2\r\,,\,\l\,)\;,
  \end{equation}
  where $\r$ is a half sum of all positive roots, and we use the scalar
  product on $\gh^*$ induced by the Killing form.

  With the Casimir element \eqref{Cas} it is associated an operator
  $\O_2$ acting on a tensor product $V'\ox V''$ of any two
  $\gg$-modules $V'$, $V''$ by the rule
  \begin{equation}\label{Cas2}
     \O_2(v'\ox v'')=\sum a_i(v')\ox b_i(v'')\;,
  \end{equation}
  which clearly does not depend on the choice of dual bases $a_i$,
  $b_i$ for $\gg$. The actions of $\O$ and $\O_2$ are related by
  the formula
  \begin{equation}\label{Cas-pol}
     \O(v'\ox v'')=\O(v')\ox v''+2\,\O_2(v'\ox v'')+v'\ox\O(v'')
  \end{equation}
  Applying this formula to $v_\hw\ox v_\hw\in V_\l\ox V_\l$, which is
  the highest vector of weight $2\l$, we get
  \begin{multline*}
    2\,\O_2(v_\hw\ox v_\hw)=(c_{2\l}-2c_\l)\cdot v_\hw\ox v_\hw=\\
    =\bigl(4(\,\l+\r\,,\,\l\,)-2\,(\,\l+2\r\,,\,\l\,)\bigr)\cdot
    v_\hw\ox v_\hw=2\,(\l,\l)\cdot v_\hw\ox v_\hw\;.
  \end{multline*}
  At the same time the formula \eqref{Cas-pol} shows that $\O_2$ does
  commute with $\gg$-action, because the Casimir element $\O$ does.
  We conclude that for any $x\in G\cdot v_\hw\subset V_\l$
  \begin{gather}
     \O(x\ox x)=c_{2\l}\cdot x\ox x=
     4\,(\,\l+\r\,,\,\l\,)\cdot x\ox x\;,
     \label{Oact}\\
     \O_2(x\ox x)=(\,\l\,,\,\l\,)\cdot x\ox x\label{O2act}\;.
  \end{gather}

\begin{proposition}\label{qe-proof-pre}
   The following four statements about $x\in V_\l$ are pairwise
   equivalent:
    \begin{enumerate}
       \item
       $x\in G\cdot v_\hw$\,;
       \item
       $x\ox x$ lays in the irreducible component
       $V_{2\l}\subset V_\l\ox V_\l$\,;
       \item
       $\O(x\ox x)=4\,(\,\l+\r\,,\,\l\,)\cdot x\ox x$\,;
       \item
       $\O_2(x\ox x)=(\,\l\,,\,\l\,)\cdot x\ox x$\,.
    \end{enumerate}
\end{proposition}

\begin{proof}
  We have seen already that  $(1)\THEN(2)\THEN(3)\IFF(4)$. The
  implication $(3)\THEN(1)$ follows from the next more precise
  statement.
\end{proof}

\begin{proposition}\label{qe-proof}
  The quadratic equations \eqref{Oact} generate the homogeneous
  ideal of the projective variety $X=\PP(G\cdot
  v_\hw)\subset\PP(V)$.
\end{proposition}

\begin{proof}
  As a $\gg$-module, the whole coordinate algebra of $\PP(V_\l)$ is
  isomorphic to $S^\bdot(V_\l^*)\simeq S^\bdot(V_\mu)$, where
  $\mu=-w_{\mathrm{max}}(\l)$ is the highest weight of $V_\l^*$ and
  $w_{\mathrm{max}}$ is the maximal length element in the Weyl
  group. It is easy to check (comp. with \cite{FH}) that its $q$-th
  homogeneous component $S^qV_\mu$ splits as
  \begin{equation}\label{Sq_dec}
     S^qV_\mu=V_{q\mu}\oplus(\O-c_{q\mu}\Id)\cdot S^qV_\mu\;,
  \end{equation}
  because the highest weight $q\mu$ appears in $S^qV_\mu$ with
  the multiplicity one and the eigenvalues of $\O$ on the irreducible
  submodules $V_\nu\subset S^qV_\mu$ with $\nu<q\mu$ are strictly
  less than $c_{q\mu}$. Let us write
  $$J_q=(\O-c_{q\mu}\Id)\cdot S^qV_\mu
  $$
  for the right summand in \eqref{Sq_dec}, which collects all
  irreducible submodules $V_\nu$ with $\nu<q\mu$. Since all the
  weights of $S^pV_\mu\cdot J_q\subset S^{p+q}V_\mu$ are
  strictly less then $(p+q)\mu$, it is clear that
  $$J=\OP_q J_q
  $$
  is a homogeneous ideal in $S^\bdot V_\mu$. We would like to check
  that $J$ is generated by $J_2$.

  The key point is that for any $v\in V_\mu$ we have
  \begin{equation}\label{OmOfvk}
     \left[\O-c_{q\mu}\Id\right](v^q)=
     \frac{q(q-1)}2\,\left[\O-c_{2\mu}\Id\right](v^2)\cdot
     v^{q-2}\;.
  \end{equation}
  Indeed, using the Leibnitz rule and the relation \eqref{Cas-pol},
  we get
  \begin{multline*}
     \O(v^q)=
        q\cdot\O(v)\cdot v^{q-1}+q(q-1)\cdot
        \O_2(v^2)\cdot v^{q-2}=\\
     =q\cdot c_\mu\cdot v^q+\frac{\,q(q-1)}2\cdot
        \(\O(v^2)-2\cdot\O(v)\cdot v\)\cdot v^{q-2}=\\
     =\(-q^2+2q\)\cdot c_\mu\cdot v^q+\frac{\,q(q-1)}2\cdot\O(v^2)\cdot v^{q-2}\;.
  \end{multline*}
  This reduces \eqref{OmOfvk} to the purely numerical identity
  $$\(q^2-2q\)c_\mu+c_{q\mu}=\frac{\,q(q-1)}2c_{2\mu}\;,
  $$
  which is verified by straightforward computation using
  \eqref{CasVal}.

  Since the powers $v^q$ span $S^qV_\mu$ as a linear space,
  the identity \eqref{OmOfvk} says actually that the ideal $J$ is
  generated by its quadratic component
  $$J_2=\(\O-c_{2\mu}\Id\)\cdot S^2V_\mu\;.
  $$
  Taking into account the equality $c_{2\l}=c_{2\mu}$, we see
  that these quadratic equations coincide with \eqref{Oact} as well
  as with ones mentioned in the condition (3) of the previous
  proposition.

  Now, rewriting the decomposition \eqref{Sq_dec} as
  $S^qV_\mu=V_{q\mu}\oplus J_q$\,, we see immediately that any
  $\gg$-invariant ideal $I\varsupsetneq J$ should contain the
  irreducible submodule $V_{q\mu}$ for all $q\GG 0$, \ie should be
  of finite codimension in $S^\bdot(V_\mu)$ as a vector space. This
  means that $J=\sqrt J$ coincides with the homogeneous ideal of the
  projective variety $Y\subset\PP(V_\l)$ defined by the quadratic
  equations \eqref{Oact}. Moreover, this means that $Y$ does not
  contain proper $\gg$-invariant closed algebraic subsets. Since
  any $G$-orbit of minimal dimension inside $Y\smallsetminus X$ would be
  such a subset, we conclude that $Y=X$.
\end{proof}

\section{Syzygies of the projective coordinate algebra}\label{par2}

\noindent
The content of this section was influenced by our discussions with
M.~Movshev. We streamline and clarify arguments used in \cite{MS2},
\cite{Mo1} and construct the Frobenius algebra structure on the
syzygies of an arbitrary projective variety $X\subset\PP(V)$
satisfying the following three properties
\begin{enumerate}
   \item
   $X$ is smooth\footnote{in fact, all results of this section (and their proofs) hold for any locally complete intersection varieties (see \cite{Ha}, \cite{Fu1} for details); all we need is well defined invertible dualizing sheaf $\o_X$};
   \item
   for $m\in\ZZ$ the cohomologies $H^i(X,\cO_X(m))=0$, if
   $i\ne 0,\,d=\dim X$\,;
   \item
   $X$ is subcanonical, \ie $\o_X=\cO_X(-N)$ for some $N\in\NN$\,.
\end{enumerate}
 As we have seen in
\S1, these conditions hold for the \SHWs
\eqref{cpe}. So, the syzygy space of any \SHW always carries a
Frobenius algebra structure.

\subsection{Coordinate algebra of a projective variety}
  By the definition, the coordinate algebra of a projective variety
  $X\subset\PP_n=\PP(V)$ is the graded algebra
  \begin{equation}\label{pca}
     A=\OP_{m\ge0}H^0(X,\cO_X(m))=S/J\;,
  \end{equation}
  where $S=\OP_{m\ge0}S^mV^*$ is the symmetric algebra of $V^*$ (the
  homogeneous coordinate algebra of $\PP(V)$) and $J=\{f\in
  S\,|\;f|_X\equiv0\}$ is the homogeneous ideal of $X$.

  In terms of generators and relations, such an algebra $A$
  is described by its {\it minimal free resolution\/}
  \begin{equation}\label{mfrA}
     \cdots\;\rTo F_2\rTo F_1\rTo F_0\rTo A\rTo0\;,
  \end{equation}
  which is an exact sequence of graded free $S$-modules of the form
  \begin{equation}\label{Rpqdef}
     F_p=\OP_{q\ge m_p}R_{p,q}\OX_\CC S[-q]\,,
  \end{equation}
  where $R_{p,q}$ are finite dimensional vector spaces of {\it $p$-th
  order syzygies of degree $q$\/} for $A$ and we write $m_p$ for the
  minimal degree appearing among the order $p$ syzygies. We will call
  $p$ the {\it homological\/} degree and $q$~--- the {\it internal\/}
  degree.

  The {\it minimality\/} of the resolution \eqref{mfrA} means that
  all homogenous components of all matrix elements of each its
  differential are polynomials of strictly positive degree, \ie for
  any $p$ the differential $F_p\rTo F_{p-1}$ takes each syzygy
  submodule $R_{p,q}\ox S[-q]$ into $\OP_{\nu\le q-1}R_{p-1,\nu}\ox
  S[-\nu]$. Thus, the tensor multiplication by the trivial $S$-module
  $\CC$ annihilates all the differentials in a minimal free
  resolution
  \eqref{mfrA} and we get for each $p$ an isomorphism of graded
  vector spaces
  \begin{equation}\label{Rpdef}
     R_p\bydef\OP_{q\ge m_p}R_{p,q}=\Tor^S_p(A,\CC)\;.
  \end{equation}
  In particular, the dimensions $\dim R_{p,q}$ do non depend on the
  choice of a minimal resolution.

\begin{example}[the HW-orbits]\label{SHWex-2}
  If $X=G/P$ is embedded into $\PP(V)$ as the highest vector orbit
  \eqref{cpe}, then its ideal $J=(Q)$ is generated by the quadratic
  equations \eqref{Oact}, which form a linear subspace
  $$Q\subset S^2V^*\subset S\;.
  $$
  Thus, the resolution \eqref{mfrA} starts with $F_0=S$, \ie $m_0=0$,
  $R_0=R_{0,0}=\CC$. Then, $F_1=Q\ox S[-2]$, \ie $m_1=2$ and
  $R_1=R_{1,2}=Q$. Further, it follows from the minimality that
  $m_p\ge p+1$ for all $p\ge1$.

  Note that the group $G$ acts naturally on the coordinate algebra
  $A$ and on the $\Tor$-spaces, thus on the syzygies. During the
  proof of the proposition \ref{qe-proof} we have seen that $Q\subset
  S^2V^*$ collects all the irreducible direct summands of $S^2V^*$
  except for $V_{2\mu}$, where $\mu$ is the highest weight of
  $V^*$. This gives an effective way to compute at least the
  starting term of minimal resolution.

  Say, for the Pl\"ucker embedding
  $\Gr(2,5)\simeq\PP(G\cdot v_{\o_2})\subset\PP(V_{\o_2})=\PP_9$\,,
  which corresponds to $\l=\o_2=[0,1,0,0]$ in the notations of \ref{g25}\,,
  we have $\mu=\o_3=[0,0,1,0]$, thus $2\,\mu=[0,0,2,0]$, and
  $S^2V^*=V_{[0,0,2,0]}\oplus V_{[0,0,0,1]}$
  (this is immediate from $\dim S^2\CC^{10}=55$). Hence,
  $Q=V_{[0,0,0,1]}$ and $\dim Q=5$.

  Similarly, for the variety of 10-dimensional pure spinors
  $\Gr^+_{\mathrm{iso}}(5,10)$ (see example \ref{so10}) we have
  $\mu=\l=\o_5$\,, $\dim V_\l=16$\,, $\dim S^2V^*=136$, but $\dim
  V_{2\mu}=126$. This implies $Q=V_{\o_1}$\,, $\dim V_{\o_1}=10$.
  Note that in the both cases, $\Gr(2,5)\subset\PP_9$ and
  $\Gr^+_{\mathrm{iso}}(5,10)\subset\PP_{15}$, the quadratic
  equations $Q$ form an irreducible $G$-module.
\end{example}

\subsection{DGA resolution for the algebra of syzygies}
The syzygies $R_{p,q}$ can be computed using the
standard Koszul resolution for the trivial $S$-module $\CC$
\begin{equation}\label{KRforC}
    \cdots\;\rTo^\dk K_2\rTo^\dk K_1\rTo^\dk K_0\rTo^\dk\CC\rTo0
\end{equation}
that has $K_p=\L^pV^*\OX_\CC S[-p]$ and the differential
$\dk=\sum\frac\PD{\PD\t_i}\ox x_i$, which takes
\begin{equation}\label{KD}
   \o\ox f\;\longmapsto\;
         \sum\frac{\PD\o}{\PD\t_i}\ox x_i\cdot f\;,
\end{equation}
where $\OFAM\t,n$ is a basis in $V^*$ considered inside the exterior
algebra $\L V^*$ and $\OFAM x,n$ is the same basis of $V^*$ but
considered inside the symmetric algebra $S V^*$. The derivation
$\PD/\PD\t_i$ takes\footnote{note that grassmann partial derivatives
skew commute
$\frac\PD{\PD\t_i}\frac\PD{\PD\t_j}=-\frac\PD{\PD\t_j}\frac\PD{\PD\t_i}$
and satisfy the graded Leibnitz rule
$\frac\PD{\PD\t_i}(\o_1\w\o_2)=\(\frac\PD{\PD\t_i}\o_1\)\w\o_2+
(-1)^{|\o_1|}\o_1\w\(\frac\PD{\PD\t_i}\o_2\)$}
$$\t_I=\IWP\t,i,k\;\longmapsto\;(-1)^{\nu-1}\t_{I\smallsetminus i}\;,\quad
  \text{when }i=i_\nu\in I\;,
$$
and annihilates all $\t_J$ with $J\not\ni i$.

Tensoring the Koszul resolution \eqref{KRforC} by $A$ over $S$, we
get a complex
\begin{equation}\label{KCforA}
    0\rTo\L^{n+1}V^*\OX_\CC A[-n-1]
     \rTo^\PD
     \;\cdots\;
     \rTo^\PD\L^2V^*\OX_\CC A[-2]
     \rTo^\PD V^*\OX_\CC A[-1]
     \rTo^\PD A
     \rTo0\;,
\end{equation}
  whose differential is given by the same formula \eqref{KD}
  considered modulo the quadratic relations $(Q)\subset S$. The
  $p$-th homology group of \eqref{KCforA} coincides with
  $$\Tor^S_p(A,\CC)=R_p
  $$
  from \eqref{Rpdef}. The Koszul complex \eqref{KCforA} can be considered from two different viewpoints. First of all, taking the direct sum of its elements, we get a DG-algebra
\begin{equation}\label{cAalgebra}
  \cA=\OP_{p\ge0}\cA_p\;,\quad\cA_p=\L^pV^*\OX_\CC A[-p]\;,
\end{equation}
whose multiplication is induced by exterior and symmetric
multiplication of the tensor factors:
$$(\o\ox f)\cdot(\eta\ox g)=(-1)^{|f||\eta|}(\o\w\eta)\ox(fg)\;,
$$
where we write $|x|\in\ZZ/2\ZZ$ for the parity of $x$ induced by
the internal degree of $x$. It satisfies $a\cdot
b=(-1)^{|a|\cdot|b|}\,b\cdot a$ for homogeneous\footnote{here and
further we write $|x|=p$ for the degree of a homogeneous element
$x\bl\cA_p$} $a\in\cA_{|a|}$, $b\in\cA_{|b|}$ and agrees with the
differential \eqref{KD}:
$$\PD(a\cdot b)=(\PD a)\cdot b+(-1)^{|a|}\,a\cdot (\PD b)\;.
$$
This implies that the syzygies $R=\OP_pR_p=H(\cA)$ also inherit an
associative algebra structure.

We are going to show that this structure is {\it graded
Frobenius\/}, \ie there is a non degenerated scalar product
$$R\OX_\CC R\rTo^{\quad a\ox b\mapsto(a,b)\quad}\CC
$$
such that $(a,b\cdot c)=(a\cdot b,c)$ and $(a,b)=\pm(b,a)$, where
the sign rule depends on the parity of $\codim X$ (see page
\pageref{sps}). This scalar product can be written as
$(a,b)=\tr(a\cdot b)$, where
\begin{equation}\label{pre-tr}
   \tr:R\rTo^{\quad a\mapsto(e,a)=(a,e)\quad}\CC
\end{equation}
is some natural {\it trace form\/} that comes from well known
geometric construction, which will be described in the next section
using another viewpoint on the Koszul complex \eqref{KCforA}.

Namely, let us fix some isomorphism\footnote{recall that $\dim
V=n+1$ and $\PP_n=\PP(V)$} $\L^{n+1}V\rTo^\d_\sim\CC$, which
performs to identify $\(\L^pV\)^*$ with $\L^{n+1-p}V$ via non
degenerated pairing
\begin{equation}\label{L-duality}
   \L^pV\ox\L^{n+1-p}V\rTo^{\quad\o\ox\eta\mapsto\o\w\eta\quad}
   \L^{n+1}V\rTo^\d_\sim\CC\;.
\end{equation}
Further, we identify $\L^p\(V^*\)$ with $\(\L^pV\)^*$ via non
degenerated pairing induced by lifting $\L^p\(V^*\)$ into
$\(V^*\)^{\ox p}$, $\L^pV$ into $V^{\ox p}$, and taking the
complete contraction\footnote{in terms of dual bases $\t_i^*\in V$,
$\t_i\in V^*$, the full contraction between $\t^*_I$ and $\t_J$
equals $\d_{IJ}\cdot\frac1{p!}$, where $p$ is the degree of the
monomials}. It is easy to see that under this identification the
Koszul complex \eqref{KCforA} turns to the complex
\begin{multline}\label{DKCforA}
    0\rTo A[-n-1]
     \rTo^\PD V\OX_\CC A[-n]
     \rTo^\PD\L^2V\OX_\CC A[-n+1]
     \rTo^\PD
     \;\cdots\\\cdots\;
     \rTo^\PD\L^{n-1}V\OX_\CC A[-2]
     \rTo^\PD\L^nV\OX_\CC A[-1]
     \rTo^\PD\L^{n+1}V\OX_\CC  A
     \rTo0\;,
\end{multline}
whose components also form an associative graded
superalgebra\footnote{in a physical cant the exchange
$\cA\leftrightsquigarrow\cA'$ is known as `odd Fourier transform'}
$$\cA'=\OP_p\L^pV\OX_\CC A[p-n-1]
$$
\wrt the multiplication induced by exterior multiplication in $\L
V$ and the usual one in $A$. The differential $\PD$ sends a
homogeneous element $a\in\cA_p$ to $p\,\Id_V\cdot a$, where
$\Id_V\in V\ox V^*$ is considered as an element of $V\ox
A[-n]\subset\cA'$ and the multiplication is taken inside $\cA'$.
Note that  $\PD$ is no longer compatible with the multiplication.
However the interpretation \eqref{DKCforA} reveals projective
geometrical meaning of $\PD$.

\subsection{Euler-\Dol bicomplex}
There is canonical Euler exact triple of coherent sheaves on
$\PP_n=\PP(V)$
$$0\rTo\cO_{\PP_n}\rTo^{\quad1\mapsto\sum\t_i^*\ox x^i\quad}
  V\ox\cO_{\PP_n}(1)\rTo T\rTo0\;,
$$
which describes the tangent sheaf $T=T_{\PP_n}$ over any point
$p\in V$ as the factor space $V/\CC\cdot p$. The exterior powers of
this triple
$$0\rTo\L^{m-1}T\rTo\L^mV\ox\cO_{\PP_n}(m)\rTo\L^mT\rTo0
$$
are naturally organized into a long exact sequence of locally free
$\cO_{\PP_n}$-modules
\begin{multline}\label{lm}
  0\rTo\cO_{\PP_n}
   \rTo V\ox\cO_{\PP_n}(1)\rTo\L^2V\ox\cO_{\PP_n}(2)\rTo\;\cdots\\
   \cdots\;\rTo\L^{n-1}V\ox\cO_{\PP_n}(n-1)\rTo\L^nV\ox\cO_{\PP_n}(n)
   \rTo\cO_{\PP_n}(n+1)\rTo0
\end{multline}
whose maps are given by the left multiplication by $\Id_V=\sum \t_i^*\ox x_i$, where
$\t_i^*\in V$ form the dual base to $\t_i\in V^*$ but $x_i\in V^*$ are considered now as global sections of $\cO(1)$, and `the multiplication' means exterior multiplication in the first factor and tensor multiplication in the second one. Hence, twisting by $\cO_{\PP_n}(k-n-1)$ and using the identification
$\L^pV\simeq\L^{n+1-p}V^*$ described above, we can rewrite \eqref{lm} as
\begin{multline}\label{lmt}
  0\rTo\cO_{\PP_n}(k-n-1)
   \rTo^\dk\L^nV^*\ox\cO_{\PP_n}(k-n)\rTo^\dk\;\cdots\\
   \cdots\;\rTo^\dk\L^2V^*\ox\cO_{\PP_n}(k-2)\rTo^\dk V^*\ox\cO_{\PP_n}(k-1)
   \rTo^\dk\cO_{\PP_n}(k)\rTo0\;,
\end{multline}
where the differential
$\L^pV^*\ox\cO_{\PP_n}(k-p)\rTo^\dk\L^{p-1}V^*\ox\cO_{\PP_n}(k-p+1)$
is given by the same formula \eqref{KD}.

Since $X$ is smooth, restricting \eqref{lmt} onto $X$
we get the following exact sequence of locally free
coherent sheaves on $X$:
\begin{multline}\label{lmtr}
  0\rTo\cO_X(k-n-1)
   \rTo^\PD\L^nV^*\ox\cO_X(k-n)\rTo^\PD\;\cdots\\
   \cdots\;\rTo^\PD\L^2V^*\ox\cO_X(k-2)\rTo^\PD V^*\ox\cO_X(k-1)
   \rTo^\PD\cO_X(k)\rTo0\;.
\end{multline}
Note that the Koszul complex \eqref{KCforA} is a direct sum of
complexes obtained from \eqref{lmtr} by applying the global
sections functor $\Gamma(X,\,\ast\,)$.

Now consider the flabby \Dol $\PDB$-resolutions\footnote{We use the
\Dol complex just by tradition. In fact, any natural flabby right
resolution $\cE^q_p$ of complex \eqref{lmtr} would be OK for the
forthcoming computation. Say, the usage of the canonical Godemant
resolution could make the computation even more transparent and works
for non smooth locally complete intersection varieties as well.} for
all coherent sheaves in \eqref{lmtr}. They are organized in the
exact bicomplex of flabby sheaves of abelian groups
$\cE^q_p=\L^pV^*\ox\O_X^{0,q}(k-p)$ on $X$:
\begin{equation}\label{EDres}
\begin{diagram}
    &&\vdots
      &&\vdots\\
    &&\uTo<{\PDB}
      &&\uTo<{\PDB}\\
    \cdots
    &\rTo^\PD
     &\L^pV^*\ox\O_X^{0,q+1}(k-p)
      &\rTo^\PD
       &\L^{p-1}V^*\ox\O_X^{0,q+1}(k-p+1)
        &\rTo^\PD
         &\cdots\\
    &&\uTo<{\PDB}
      &&\uTo<{\PDB}\\
    \cdots
    &\rTo^\PD
     &\L^pV^*\ox\O_X^{0,q}(k-p)
      &\rTo^\PD
       &\L^{p-1}V^*\ox\O_X^{0,q}(k-p+1)
        &\rTo^\PD
         &\cdots\\
    &&\uTo<{\PDB}
      &&\uTo<{\PDB}\\
    &&\vdots
      &&\vdots
\end{diagram}
\end{equation}
Writing $\cE^q_p$ in the $(-p,q)$-cell of the second coordinate
quadrant of $(p,q)$-plane, we get a diagram bounded by inequalities
$-n-1\le-p\le0$, $0\le q\le d$, where $d=\dim X$. It has the
following obvious properties:
\begin{enumerate}
  \item
   All the rows of \eqref{EDres} are exact and acyclic \wrt the
   functor $\Gamma(X,\,\ast\,)$. In particular, taking the global
   sections in \eqref{EDres}, we get a bicomplex
   $\Gamma(X,\cE^q_p)$, with exact rows, whose associated total
   complex is exact.
  \item
  For each $p$ the $p$-th column of \eqref{EDres} gives a flabby
  resolution for the coherent sheaf $\L^pV^*\ox\cO(k-p)$. So,
  applying $\Gamma(X,\,\ast\,)$ to $p$-th column of \eqref{EDres},
  we get a complex whose $q$-th cohomology group equals $\L^pV^*\ox
  H^q(X,\cO_X(k-p))$.
\end{enumerate}
  Hence, there is a spectral sequence that converges to the zero cohomologies of the total complex associated with $\Gamma(X,\cE^q_p)$ and it has
\begin{equation}\label{E1ED}
   E_1^{-p,q}=H^q\bigl(\Gamma(X,\cE^q_p)\,,\,\ovl\PD\,\bigr)
             =\L^pV^*\ox H^q(X,\cO_X(k-p))\;.
\end{equation}
  Now, if $X$ satisfies the condition (2) formulated on page
  \pageref{par2}, then
  $$H^q(X,\cO_X(k-p))=0\quad\text{for }q\ne0,d\;.
  $$
  Therefore all non zero terms of the spectral sequence \eqref{E1ED} will be situated only in two horizontal rows: $q=0$ and $q=d=\dim X$.

\begin{proposition}\label{EdLem}
   Under the assumptions\/ {\upshape(1\,--\,3)} on page \pageref{par2} for each
   $k\in\ZZ$ and any $0\le p\le n+1$ there is an isomorphism
   \begin{equation}\label{duadif}
      \tau_p\,:\;R_{n-d-p,n-N+1-k}^*\rTo^\sim R_{p,k}
   \end{equation}
   provided by the differential in $E_d$-term of spectral sequence
   \eqref{E1ED}.
\end{proposition}

\begin{proof}
  The bottom row $q=0$, of \eqref{E1ED}, coincides with the internal degree $k$ homogeneous slice of the Koszul resolution \eqref{KCforA}:
  $$0\rTo\L^{n+1}V^*\OX_\CC A_{k-n-1}\;
       \rTo^\PD\cdots
       \rTo^\PD\L^2 V^*\OX_\CC A_{k-2}
       \rTo^\PD V^*\OX_\CC A_{k-1}
       \rTo^\PD A_k
       \rTo0\;.
  $$
  Hence, the bottom row in the $E_2$-term consists of the following
  syzygies:
$$E_2^{-p,0}=\Tor_S^p(A,\CC)_{k}=R_{p,k}\;.
$$
For the upper row $q=d$ in \eqref{E1ED} we claim that it is dual to
the complex
\begin{multline*}
    0\lTo A_{n-N+1-k}
     \lTo^\PD V^*\OX_\CC A_{n-N-k}
     \lTo^\PD\L^2 V^*\OX_\CC A_{n-N-k-1}\lTo^\PD\;\cdots\\
     \cdots\;\lTo^\PD\L^{n+1}V^*\OX_\CC A_{-N-k}
     \lTo0\;,
\end{multline*}
because of $H^d(X,\cO_X(k-p))\simeq H^0(X,\cO_X(p-k-N))^*$ by the
Serre duality\footnote{here we use the condition (3) on page
\pageref{par2}\:: $\o_X=\cO_X(-N)$} and
$\L^pV^*\simeq\(\L^{n+1-p}V^*\)^*$ via the dual version of pairing
\eqref{L-duality}. So, the top row in $E_2$ is filled by the spaces
dual to the following syzygies:
$$E_2^{-p,d}=\Tor_S^{n+1-p}(A,\CC)_{n-N+1-k}^*=R_{n+1-p,n-N+1-k}^*\;.
$$
The differentials in the consequent terms of this spectral sequence
will be non trivial only in $E_d$. Therefore, to get the zero limit, the
$E_d$-differential should map $(-p-d-1,d)$-cell isomorphically onto
$(-p,0)$-cell providing isomorphism
\eqref{duadif}.
\end{proof}

\begin{corollary}
   If the ideal of $X$ is generated by some linear space of
   quadrics $Q\subset S^2(V^*)$, then non zero syzygies $R_{p,k}$
   can appear only for the following values $p$, $k$:
   \begin{itemize}
     \item
     $p=0$, $k=0$ or $p=n-d$, $k=n-N+1$, where
     $R_{0,0}=\CC\simeq R_{n-d,n-N+1}^*$\,;
     \item
     $p=1$, $k=2$ or $p=n-d-1$, $k=n-N-1$, where
     $$R_{1,2}=Q\simeq R_{n-d-1,n-N-1}^*\;;
     $$
     \item
     $2\le p\le n-d-2$ and $p+1\le k\le p+d-N$, where
     $R_{p,k}\simeq R_{n-d-p,n-N+1-k}^*$\,.
   \end{itemize}
\end{corollary}

\begin{proof}
  Indeed the syzygies $R_{p,k}=0$ automatically vanish in the
  following four cases:
  $$\begin{aligned}
       1)&\quad p<0\quad\text{or}\quad p>n+1
       &\qquad2)&\quad p=0\quad\text{and}\quad k\ne0\\
       3)&\quad p=1\quad\text{and}\quad k\ne2
       &\qquad4)&\quad 2\le p\le n+1\quad\text{and}\quad k\le p
    \end{aligned}
  $$
  Applying this to the right hand side of \eqref{duadif} we see
  that $R_{p,k}=0$ also for $p>n-d$ or $p<d-1$, for $p=n-d$ and
  $k\ne n-N+1$, for $p=n-d-1$ and $k\ne n-N-1$, and finally, for
  $2\le p\le n-d-2$ and $k\ge p+d-N-1$.
\end{proof}

\begin{example}[continuation of \ref{g25}, \ref{so10}, and \ref{SHWex-2}]\label{SHWex-2c}
  For the grassmannian $\Gr(2,5)\subset\PP_9$ we have $n=9$, $d=6$, $N=5$. Thus,
  by the above corollary,
  $$R_{3,5}^*\simeq R_{0,0}=\CC\;,\qquad
    R_{2,3}^*\simeq R_{1,2}=Q\simeq\CC^5
  $$
  and all the other syzygies vanish. For 10-dimensional pure spinors
  $\Gr^+_{\mathrm{iso}}(5,10)\subset\PP_{15}$ we have $n=15$, $d=10$, $N=8$. The previous corollary
  implies that
  $$R_{5,8}^*\simeq R_{0,0}=\CC\;,\qquad
    R_{4,6}^*\simeq R_{1,2}=Q\simeq\CC^{10}
  $$
  and all other syzygies of orders $0,\,1,\,4,\,5$ vanish. The
  precise computation of $R_2\simeq R_3^*$ requires more
  sophisticated computational analysis of the
  $\mathrm{Spin}(10,\CC)$-module $R$. It was made, e.\;g. in
  \cite{CN}, \cite{CR}, \cite{MS2}.
\end{example}

\subsection{\mathversion{bold}Scalar product and trace on $R$}\label{sps}
Let us define the trace functional \eqref{pre-tr} as a linear form
on $R$ that annihilates all $R_{p,k}$ except for $R_{n-d,n-N+1}$
and sends $R_{n-d,n-N+1}$ to $R_{0,0}^*=\CC$ via isomorphisms $\tau_0^{-1}$ inverse to $\tau_0$ defined in \ref{EdLem}.
This form provides the syzygy algebra $R$ with a scalar product
\begin{equation}\label{SPonR}
  (a,b)\bydef\tr(a\cdot b)
\end{equation}
which satisfies the property $(a\cdot b,c)=(a,b\cdot c)=\tr(a\cdot
b\cdot c)$, because $R$ is associative.

Since for $a\in R_{p,k}$, $b\in R_{n-d-p,n-N+1-k}$ we have $a\cdot b=(-1)^{p(n-d-p)}b\cdot a$\,, the scalar product \eqref{SPonR} is purely symmetric, if $\codim_{\PP_n}X=(n-d)$ is odd. If $\codim_{\PP_n}X$ is even, then \eqref{SPonR} is {\it even supersymmetric\/}, that is symmetric, when the both arguments are even, and skew-symmetric, when the both arguments are odd.

\begin{proposition}
  The scalar product \eqref{SPonR} is non degenerate, and $\fain aR$
  $$\tau_0^{-1}(a)=(a,\ast)
  $$
  as the linear forms on $R$.
\end{proposition}

\begin{proof}
There is a natural scalar product on $E_1$-term of the spectral
sequence \eqref{E1ED} induced by the multiplication and the trace
form provided by the Serre duality. It takes
\begin{equation}\label{SPonE1}
\begin{diagram}
  \L^{p_1}V^*\ox H^0\(\cO_X(m_1)\)\x\L^{p_2}V^*\ox H^d\(\cO_X(m_2)\)\\
  \dTo\\
  \L^{n+1}V^*\ox H^d\(\cO_X(-N)\)\text{\rlap{$\simeq\CC$}}
\end{diagram}
\end{equation}
for $p_1+p_2=(n+1)$ and $m_1+m_2=-N$ (all the other components of
$E_1$ are mutually orthogonal). Clearly, it is non degenerated. We
extend the isomorphism from the bottom row of \eqref{SPonE1} to a linear map $E_1\rTo^\Tr\CC$\,,
which annihilates all the components of $E_1$
except $\L^{n+1}V^*\ox H^d\(\cO_X(-N)\)$, and define the scalar
product on $E_1$ by prescription
\begin{equation}\label{E1pairing}
   (a,b)_1\bydef\Tr(ab)\;.
\end{equation}
Since $d_1$ satisfies $d_1(ab)=(d_1a)\,b+(-1)^{|a|}\,a\,(d_1b)$ and
$\im d_1$ is annihilated by $\Tr$, this scalar product interacts
with $d_1$ by the rule
\begin{equation}\label{SPonE1sp}
   \(d_1a\,,\,b\)_1=\Tr\((d_1a)\,b\)=
   (-1)^{|a|+1}\Tr\(a\,(d_1b)\)=
   (-1)^{|a|+1}\(a\,,\,d_1b\)_1\;.
\end{equation}

Let us show that this scalar product induces well defined non
degenerate pairing on the cohomologies $H(E_1,d_1)=E_2=R$ and that
this pairing coincides with \eqref{SPonR} and produces the duality
maps \eqref{duadif}. We can fix some vector space decomposition
(compatible with the $p$-grading on $E_1$)\,:
$E_1=Z\oplus W=d_1W\oplus C\oplus W$
such that $d_1$ takes $W$ isomorphically onto $d_1W=\im(d_1)$,
$Z=d_1W\oplus  C=\ker(d_1)$, and $C\simeq E_2$ consists of
representatives for the cohomology classes. It follows from
\eqref{SPonE1sp} that $Z\subset\(d_1W\)^\bot$. Hence, the scalar
product of cohomology classes is well defined, and for each $p$ the
pairing
$$\frac{W_p}{W_p\cap\(d_1W\)^\bot}\x d_1W_{n-p}\rTo\CC
$$
induced by \eqref{E1pairing} is non degenerate\footnote{note that
$d_1W_{n-p}\subset\(d_1W\)_{n+1-p}$}. Let
$$\begin{aligned}
     w_p&=\dim_\CC W_p=\dim_\CC d_1W_p\,,\\
     w'_p&=\dim_\CC\Bigl(W_p/\bigl(W_p\cap(d_1W)^\bot\bigr)\Bigr)\,.
  \end{aligned}
$$
Then $w'_p=w_{n-p}$ for each $p$ and evident inequalities $w'_p\le
w_p=w'_{n-p}\le w_{n-p}=w'_p$ imply that $w'_p=w_p$ for all $p$.
So, $W\cap\(d_1W\)^\bot=0$, that is $\(d_1W\)^\bot=Z$ and
$Z^\bot=\(d_1W\)$. This implies $C^\bot\cap C=0$, which means that
\eqref{E1pairing} gives a non degenerate pairing on $C$.

Now write $\o\in \L^{n+1}V^*\ox H^d\(\cO_X(-N)\)$ for the basic
element being sent to 1 by  the isomorphism from the bottom row of \eqref{SPonE1}.
For any $a\in E_2^{-p,0}$ there exists some $b\in E_2^{-p-d-1,d}$
such that $\Tr(ab)=1$. Then $ab=\o$ in $R_{n-d,n-N+1}$.
Since in $E_d$-term       $d_{E_d}(a)=0$, we have
$1=d_{E_d}(\o)=d_{E_d}(ab)=(-1)^{p}\,a\,d_{E_d}(b)$\,.
Therefore $\tr(a\,d_{E_d}(b))=1$, which means that for any $a\in R$
the linear form $\tr(a\cdot *)$ is non zero and coincides with
$\tau_0^{-1}(a)$.
\end{proof}

\section{Cohomology of the dual graded Lie superalgebra}

  We are going to compare an algebra of syzygies for an arbitrary commutative graded quadratic Koszul algebra $A$ and an algebra of cohomologies of a graded Lie superalgebra
  $L$ Koszul dual to $A$ in the sense of Ginzburg and Kapranov
  \cite{GK}. Namely, in \ref{HcompT} we identify the syzygies of
  $A$ with the cohomologies of Lie subalgebra $\Lt\subset L$ and
  give alternative description for the algebra structure on the
  syzygies.

  By quite deep theorem of R.~Bezrukavnikov (see \cite{Bz}) the
  projective coordinate algebra of any highest weight orbit $X=G/P$
  (not necessary subcanonical) is Koszul. Thus our results can be
  applied to the syzygies of the highest weight vector orbits.

  Certainly, the coincidence of two algebra structures on the space of
  syzygies (one constructed in \S2 and another we will construct in
  this section) could be extracted from the general bar-cobar
  equivalence staff\footnote{see \cite{Ke1}, \cite[sec.\;4]{Ke2}
  for the most comprehensive approach}. But in our situation this
  will be clearly apparent, fortunately, and we will check it `by
  hands'.  We begin with recallment of some standard resolutions
  and the experts could jump directly to \ref{CH-summary}.

\subsection{Dual quadratic algebra}\label{QDA}
  Recall that a {\it quadratic algebra\/}\footnote{recall that in
  this paper we restrict ourself by $\CC$-algebras only; but in
  this section the reader can everywhere replace $\CC$ by an
  arbitrary field of zero characteristic} generated by a vector
  space $V^*$ is an associative algebra $A$ of the form
  $$A=\TT(V^*)/(I)\;,
  $$
  where $\TT(V^*)$ is the tensor algebra of a vector space $V^*$ and
  $(I)\subset\TT(V^*)$ is a double side ideal spanned by a vector
  subspace $I\subset V^*\ox V^*$, of homogeneous quadratic relations.
  Well known Priddy's construction (see \cite[p.\;108]{GM}, \cite{Pd})
  attaches to any such an algebra $A$ the dual quadratic algebra
  $A^!$ generated by the dual space $V$ with the orthogonal relation
  ideal
  $$A^!=\TT(V)/(I^\bot)\;,
  $$
  where $I^\bot\subset V\ox V$ is the annihilator of $I$. Clearly,
  $A^{!!}=A$.

  A projective coordinate algebra $A=\OP_{m\ge0}H^0(X,\cO(m))$ of
  any variety $X\subset\PP(V)$ whose ideal is generated by
  quadratic equations $\{q^\nu\}\subset S^2V^*$ fits into this
  framework as
  $$A=S^\bdot(V^*)/(Q)=\TT(V^*)/(C+Q)\;,
  $$
  where $C=\mathrm{Skew}(V^*\ox V^*)\simeq\L^2V^*$ consists of
  commutativity relations and $Q\subset\Sym(V^*\ox V^*)$ is the
  linear span of symmetric bilinear forms $\wtd{q^\nu}$, polarizing
  the quadratic equations for $X$. In this case the Priddy dual
  algebra
  $$A^!=\TT(V)/(Q^\bot\cap\Sym(V\ox V))
  $$
  can be tautologically treated as universal enveloping algebra
  for graded Lie superalgebra\footnote{this is the simplest
  motivating example for much more wide Koszul duality between
  graded Com and Lie operads, see \cite{GK}}
  \begin{equation}\label{Ldef}
     L=\OP_{m\ge 1} L_m=\Lie(V)/(\Ann(Q))\;,
  \end{equation}
  which is a factor of free graded Lie s-algebra generated by
  $V$ (taken with {\it odd\/} parity) through graded Lie ideal
  generated by $\Ann(Q)\subset\Sym(V\ox V)$\,, \ie by
  $$\left\{\,v_1\ox v_2+v_2\ox v_1\,|\;
           \wtd q(v_1,v_2)=0\;\fain qQ\right\}\;.
  $$
  Thus, $L_1=V$\,,\quad$L_2=S^2V/\Ann(Q)\simeq Q^*$\,, \etc

  In terms of coordinates, if we fix some dual bases $\{v_i\}$,
  $\{x^i\}$ for $V$, $V^*$ and $\{q^\nu\}$, $\{z_\nu\}$ for $Q$,
  $Q^*$, then we can describe $L$ as graded Lie s-algebra
  generated by $L_1$ with elements $v_i$ of parity 1 as a basis for
  $L_1$, elements $z_\nu$ as a basis for $L_2$, and Lie s-brackets
  given by
  \begin{equation}\label{Lbrac}
    [v_i,v_j]=\sum_\nu\wtd{q^\nu}(v_i,v_j)\cdot z_\nu=
    \sum_\nu a^\nu_{ij}z_\nu\;,
  \end{equation}
  where $a^\nu_{ij}$ is the matrix of $q^\nu$ in the basis $v_i$,
  \ie $q^\nu=\sum a^\nu_{ij}x^ix^j$.

  \subsection{Bar construction}
  Recall that for any graded associative $\CC$-algebra
  $B$ with unity and augumentation $B\rOnto^\e\CC$, there
  is the bar-complex of free graded left $B$-modules
  \begin{equation}\label{bar}
    \cdots\;\rTo^{\ovl\PD}
    B\ox[B^{\ox3}]\rTo^{\ovl\PD}
    B\ox[B^{\ox2}]\rTo^{\ovl\PD}
    B\ox[B]\rTo^{\ovl\PD}
    B\rTo^\e\CC\rTo0\;,
  \end{equation}
  where the tensor products are taken over $\CC$ and the $B$-linear
  differential is defined on the free generators by prescription
  \begin{multline}\label{dbar}
     \ovl\PD(1\ox[b_1\ox b_2\ox\dots\ox b_m])=
     b_1\ox[b_2\ox\dots\ox b_m]+\\
     +1\ox\sum_{i=1}^{m-1}(-1)^i[b_1\ox\dots\ox(b_ib_{i+1})
                             \ox\dots\ox b_m]\;.
  \end{multline}
  It is clearly contracted onto $\CC$ by the homotopy taking
  \begin{equation}\label{hbar}
      b_0\ox[b_1\ox\dots\ox b_m]\mapsto
      1\ox[b_0\ox b_1\ox\dots\ox b_m]
  \end{equation}
  and gives the standard free resolution for the trivial $B$-module
  $\CC$ in the category of graded left $B$-modules. The space of free
  generators $\TT^c(B)$ carries the natural coalgebra structure
  dual to the tensor multiplication
  $$\begin{gathered}
    \TT^c(B)\rTo^{\quad\D\quad}
    \TT^c(B)\ox\TT^c(B)\\
    [b_1\ox\dots\ox b_m]\mapsto
    \sum_{i=0}^{m}[b_1\ox\dots\ox b_i]\,\ox\,
        [b_{i+1}\ox\dots\ox b_m]
  \end{gathered}
  $$
  (where $[]\bydef1$) and the bar differential is a coderivation
  \wrt this coproduct, \ie satisfies\footnote{applying a homogeneous operator monomial
  $f_1\ox f_2\ox\dots\ox f_m$ to a homogeneous vector monomial
  $v_1\ox v_2\ox\dots\ox v_m$ we always assume the Koszul
  sign agreements:
  $f_1\ox f_2\ox\dots\ox f_m(v_1\ox v_2\ox\dots\ox v_m)=
    (-1)^\e\cdot f_1(v_1)\ox f_2(v_2)\ox\dots\ox f_m(v_m)$\,,
  where $\e=|f_m|\cdot(|v_1|+\dots+|v_{m-1}|)+
            |f_{m-1}|\cdot(|v_1|+\dots+|v_{m-2}|)+\dots+
            |f_2|\cdot|v_1|$}
  $$(1\ox\ovl\PD+\ovl\PD\ox1)\comp\D=\D\comp\ovl\PD\;.
  $$
  Thus, $\ext_B^\bdot(\CC,\CC)$ can be described as the cohomology
  algebra of the DG algebra
  \begin{equation}\label{cobar}
     \hom_B(B\ox\TT^c(B),\CC)=\hom_\CC(\TT^c(B),\CC)
     =\TT(B^*)
  \end{equation}
  whose multiplication is the standard tensor
  multiplication\footnote{it is instructive to see how does it
  agree with the classic Yoneda product
  $$\ext^k_B(\CC,\CC)\ox\ext^m_B(\CC,\CC)
    \rTo^{\quad\f\ox\psi\mapsto\f\comp\psi\quad}
    \ext^{k+m}_B(\CC,\CC)
  $$
  defined by obvious extending of $\psi\in{B^{\ox m}}^*$,
  $\f\in{B^{\ox m}}^*$ to the $B$-linear homomorphisms
  $$\wtd\f\,,\,\wtd\psi\,:\;B\ox\TT^c(B)\rTo\CC\;,
  $$
  then lifting $\wtd\psi$
  to some degree $m$ homomorphism of free resolutions
  $B\ox\TT^c(B)\rTo^{\wtd\psi_\bdot}B\ox\TT^c(B)$, and taking the
  composition $1\ox B^{\ox(m+k)}\rTo^{\;\wtd\psi_{m+k}\;}
  B\ox B^{\ox k}\rTo^{\e\ox1}1\ox B^{\ox k}\rTo^\f\CC$\;;
  clearly, $\wtd\psi_\bdot=(1\ox\psi)\comp\D$ gives precisely the
  required lifting of $\psi$ to the morphism of bar resolutions and
  leads to the Yoneda product $\psi\comp\f=(\psi\ox\f)\comp\D$,
  which coincides with the tensor product of multilinear forms} and
  the differential is dual to \eqref{dbar}, \ie takes 1 to zero,
  acts on degree 1 generators $\b\in B^*$ as
  $$\ovl\PD^*\b(b_1,b_2)=\b(b_1b_2)\;,
  $$
  and is extended onto the whole of $\TT(B^*)$ by the Leibnitz
  rule. We call \eqref{cobar} the {\it cobar complex\/} of $B$. It
  is naturally bigraded. In what follows we always call the degree
  \wrt the natural grading in the tensor algebra as {\it the
  (co)\,homological degree\/} in a contrast with {\it the internal
  degree\/}, which equals the total sum of degrees of all tensor
  factors \wrt the internal grading of $B$.

  \subsection{Koszulity.}\label{KozCrit}
  Let $B=A^!=\TT(V)/(I^\bot)$ be the dual quadratic algebra for
  $A=\TT(V^*)/(I)$. There is the Koszul complex of graded
  left $B$-modules
  \begin{equation}
     K_B=(B\ox A^*,\dk)
  \end{equation}
  whose differential $\dk$ comes from the right $B\ox A$-module\footnote{the algebra structure on $B\ox A$ is given by $(a_1\ox b_1)\cdot(a_2\ox b_2)=(-1)^{|a_2||b_1|}(a_1a_2)\ox(b_1b_2)$, where $|x|$ means the internal degree of $x$ modulo 2} structure on $B\ox A^*$ given by the right multiplication in $B$ and dual to the left multiplication in $A$. Namely, it is easy to see that the Casimir element
  $$\Id_V\in\End_\CC(V)=V\ox V^*=B_1\ox A_1\subset B\ox A
  $$
  has the zero square in $B\ox A$. By the definition, $\dk$ is given by the right
  action of $\Id_V$ on $B\ox A^*$. In `low level'
  notations, if $v_i$, $x_i$ are dual bases for $V$ and $V^*$, then
  $$\dk(b\ox\a)=\sum (b\cdot v_i)\ox(\a\comp x_i)\;,
  $$
  where\quad
  $\a\comp x_i=\(A\rTo^{\;a\mapsto\a(x_i\cdot a)\;}\CC\)\in A^*$\,.
  For example, if
  \begin{gather*}
    B=\L(V)=\TT(V)/\Sym(V\ox V)\;,\\
    A=S(V^*)=\TT(V^*)/\mathrm{Skew}(V^*\ox V^*)
  \end{gather*}
  are the ordinary exterior and symmetric algebras, then
  $$K_{S(V)}=\(\,S(V)\ox\L(V)\,,\;\sum v_i\ox\frac\PD{\PD v_i}\,\)
  $$
  is the Koszul complex \eqref{KRforC} but with $V$ instead of
  $V^*$.

  There is canonical morphism of the differential graded $B$-modules
  \begin{equation}\label{kos-bar}
     K_B=B\ox A^*\rTo^{\quad1\ox\tau\quad}B\ox\TT^c(B)
  \end{equation}
  induced by the coalgebra morphism
  $A^*\rTo^\tau\TT^c(V)=\TT^c(B_1)\subset\TT^c(B)$ dual to
  the structure morphism of algebras $\TT(V^*)\rOnto
  A=\TT(V^*)/(I)$. It is well known and not difficult to check (see
  \cite{PP}, \cite{Pd}) that the following conditions on $B$ are
  pairwise equivalent:
\begin{enumerate}
   \item
   the Koszul complex $K_B$ gives a free graded left $B$-module
   resolution for $\CC$\,, \ie the mapping \eqref{kos-bar} is a
   quasiisomorphism;
   \item
   $A\simeq\ext^\bdot_B(\CC,\CC)$\,;
   \item
   $\ext^{i,j}_B(\CC,\CC)=0$ for $i\ne j$, where $\ext^{i,j}$ means
   the internal degree $j$ graded component of $i$-th derived
   functor $\ext^i$\,;
   \item
   for each $m\ge3$ the subspaces $W_\nu=V^{\ox \nu}\ox I^\bot\ox V^{\ox(m-\nu-2)}\subset V^{\ox m}$ (where $0\le\nu\le(m-2)$) form a distributive lattice\footnote{recall that this means the coincidences $W_\a\cap(W_\b+W_\g)=W_\a\cap W_\b+W_\a\cap W_\g$ and $W_\a+(W_\b\cap W_\g)=(W_\a+W_\b)\cap(W_\a+W_\g)$ for all $\a$, $\b$, $\g$ or, equivalently, the existence of a basis $E=\{e_i\}\subset V^{\ox m}$ for $V^{\ox m}$ such that $\forall\a$ $W_\a\cap E$ is a basis for $W_\a$} in $V^{\ox m}$\,.
\end{enumerate}
  Quadratic algebras satisfying this conditions are called {\it
  Koszul algebras\/}. Of course, $B$ is Koszul iff $A=B^!$ is
  Koszul, and one can exchange $B$ and $A$ in the above properties.
  Since $\hom_B(\,\ast\,,\,\CC)$ kills the Koszul differential,
  applying this functor to \eqref{kos-bar} we get DGA homomorphism
  from $A$, considered as DG algebra with the zero differential, to
  the cobar complex \eqref{cobar}
  \begin{equation}\label{A-quism}
     A\rTo^{\tau^*}(\TT(B^*),\ovl \PD^*)\;,
  \end{equation}
  which is a quasiisomorphism as soon as $A$, $B$ are Koszul. Thus,
  each Koszul algebra $A$ is {\it canonically\/} identified with
  the algebra $\ext_B(\CC,\CC)$ via \eqref{A-quism}.

\subsection{Chevalley complex.}
  For a commutative Koszul algebra $A=S(V^*)/(Q)$ the quaisiisomorphism \eqref{A-quism} means the natural identification of $A$ with the Lie algebra cohomology
  $$A\simeq\ext^\bdot_{U(L)}(\CC,\CC)\bydef H^\bdot(L,\CC)\;,
  $$
  which can be computed using another reduction of the bar complex for $B=U(L)$ known as Chevalley's complex. Let us write $\L^c(L)$ for
  the graded s-exterior coalgebra, which is dual to the tensor
  algebra of $L^*$ factorized through the relations $u\ox
  v+(-1)^{|u||v|}v\ox u$. It can be considered as the sub-coalgebra
  of $\TT^c(B)$ via the s-alternation embedding
  \begin{equation}\label{s-alt}
    e_1\w e_2\w\dots\w e_m\mapsto
    \frac1{m!}\sum_{\s\in\gS_n}\text{s-}\mathrm{sgn}(\s)\,
    e_{\s(1)}\w e_{\s(2)}\w\dots\w e_{\s(m)}
  \end{equation}
  where s-sign takes proper account of the internal degree of the
  permuted elements. One can check (see \cite[\S3\;ex.21]{Bu1},
  \cite[ch.XIII,\;ex.14]{CE}) that the bar differential on $\TT^c(B)$
  takes the subcoalgebra $\L^c(L)$ to itself and the resulting
  subcomplex\footnote{whose differential $\dc$ is the restricted
  bar differential}
  \begin{equation}\label{ChC-display}
    \cdots\;\rTo^{\dc}
    B\ox\L^3 L\rTo^{\dc}
    B\ox\L^2 L\rTo^{\dc}
    B\ox L\rTo^{\dc}
    B\rTo^\e\CC\rTo0
  \end{equation}
  gives the free graded left $B$ module resolution for $\CC$ as
  well. We call it {\it the Chevalley resolution\/} and denote by
  $\gC$ or $\gC(L)$ when the precise reference to $L$ is important.

  Practical handling of \eqref{ChC-display} becomes more
  demonstrative with an alternative Lie theoretic interpretation of
  $\gC$. Namely, let us write $\bar L$ for another copy of the
  vector superspace $L$ but with the inverse parity and the trivial
  abelian s-Lie structure. Then we can write
  $\gC_m=\L^m(L)=S^m(\ovl L)$ (because of the parity change) and
  treat
  $$\gC=B\ox\L^m(L)=B\ox S^m(\ovl L)=U(L\oplus\Lb)
  $$
  as the universal enveloping algebra of an abelian extension
  $L\oplus\ovl L$ that contains $L$ as a Lie subalgebra,
  $\Lb$ as an abelian ideal, and has Lie brackets defined by
  prescriptions\footnote{of course, the last two formulas,
  describing the action of $L$ on $\Lb$, are equivalent, because of
  $[\ovl y,x]=-(-1)^{|x|(1+|y|)}[x,\ovl y]=
   -(-1)^{|x|(1+|y|)}\ovl{[x,y]}=
   (-1)^{|x|(1+|y|)}(-1)^{|x||y|}\ovl{[y,x]}$}
  \begin{equation}\label{LLrel}
  \begin{aligned}
     {[x,y]}&=[x,y]_ L\;, &\quad [\ovl x,\ovl y]&=0\;,\\
     [x,\ovl y]&=\ovl{[x,y]}_ L\;,
              &\quad [\ovl y,x]&=(-1)^{|x|}\ovl{[y,x]}_ L
  \end{aligned}
  \end{equation}
  for all $x,y\in L$. Thus, in new s-symmetric notations, $\gC_m$
  consists of degree $m$ s-symmetric monomials
  \begin{equation}\label{Lbmons}
     b\cdot \ovl e_1\ovl e_2\,\cdots\,\ovl e_m\;,\quad
     \text{where }b\in B=U(L)=\gC_0\,,\quad e_i\in L\;.
  \end{equation}
  They also can be embedded into $\TT^c(L)$ via s-symmetrization
  identical with \eqref{s-alt}.

  To get an alternative description for the Chevalley differential
  \eqref{ChC-display}, let us define for a moment a new
  differential\footnote{we will see soon that it coincides with
  $\dc$} $\gC\rTo^d\gC$ as the odd right s-algebra
  derivation\footnote{\ie satisfying the right Leibnitz rule
  $d(ab)=ad(b)+(-1)^{|b|}d(a)b$} whose action on the generating
  vector space $L\oplus\Lb$ is given by by the same rule as the bar
  differential
  \begin{equation}\label{ChDdef}
     d(x)=0\;,\qquad d(\ovl x)=x\qquad\fain xL\;.
  \end{equation}
  It is clear that
  $d$ preserves the enveloping algebra relations and
  automatically satisfies the right Leibnitz rule \wrt
  supercommutators:
  $$d([a,b])=[a,d(b)]+(-1)^{|b|}[d(a),b]\;.
  $$
  Its action on the generators \eqref{Lbmons} looks like
  \begin{equation}\label{ChD_in_coords}
     d\(\ovl e_1\ovl e_2\,\cdots\,\ovl e_m\)=
        \sum_{1\le j\le m}\!\!\!\pm\,e_j\cdot\ovl e_1\,\cdots\,
        \without{\ovl e}j\,\cdots\,\ovl e_m+\!\!\!
        \sum_{1\le i<j\le m}\!\!\!\pm\,\ovl{[e_i,e_j]}\,
        \ovl e_1\,\cdots\,\without{\ovl e}i\,\cdots\,
                          \without{\ovl e}j\,\cdots\,
        \ovl e_m\;,
  \end{equation}
  where the precise sign calculation is quite cumbersome, but it is
  not so important for our purposes\footnote{for example, to see
  that $d^2=0$, it is enough to mention that $d$, being an odd
  right derivation, forces the commutator $[d,d]=2\,d^2$ to be the
  right derivation of $\gC$ as well; now $d^2=0$ follows from
  \eqref{ChDdef}}. Note that the coalgebra structure on the bar
  complex agrees with the standard coalgebra structure on the
  universal enveloping algebra, which is given on generators
  $\ell\in L\oplus\Lb$ by the usual rule
  \begin{equation}\label{DeltaDef}
     \D(\ell)=1\ox\ell+\ell\ox1\;,
  \end{equation}
  and is extended onto the whole of $\gC$ as a homomorphism of
  graded algebras
  $$\gC\rTo^\D\gC\ox\gC\,,
  $$
  where the algebra
  structure $\gC\OX_\CC\gC$ is given by
  $(a\ox b)\cdot(c\ox d)=
    (-1)^{|b||c|}(a\cdot c)\ox(b\cdot d)$\,.
  The differential $d$
  agrees with the coalgebra structure \eqref{DeltaDef}, \ie
  satisfies\footnote{indeed, since $\D$ is an algebra homomorphism
  and both $d$, $(1\ox d+d\ox1)$ are the right derivations,
  the both sides of \eqref{coLeib} are right
  derivations of $\gC$ with values in $\gC\ox\gC$; by
  \eqref{ChDdef} and \eqref{DeltaDef} they coincide on the
  generating vector space $L\oplus\Lb$}
  \begin{equation}\label{coLeib}
     \D\comp d=(1\ox d+d\ox1)\comp\D\;.
  \end{equation}
  We conclude that $d=\dc$ on $\gC\subset\TT(B)$. In particular,
  this gives the `low level' description for
  \eqref{ChC-display} via \eqref{ChD_in_coords}.

  Thus, we can
  compute $\ext_B(\CC,\CC)$ as the cohomologies of the complex
  $$\(C^\bdot(L),d\)=\hom_B(\,\gC(L)\,,\,\CC)
                      =\hom_\CC(\L(L),\CC)=\(\L(L^*),\dc^*\)\;,
  $$
  which carries the natural DG algebra structure whose
  multiplication is induced by the multiplication in the s-exterior
  algebra and the differential is induced by \eqref{ChD_in_coords}.
  Let us finalize this preliminary discussion as

\begin{proposition}\label{CH-summary}
  Let $A$ be an arbitrary commutative Koszul quadratic algebra,
  $B=A^!=U(L)$ be its dual algebra treated as the universal
  enveloping algebra for a graded Lie s-algebra $L$. Then $A$,
  considered as DG algebra with the zero differential, admits
  canonical isomorphism
  \begin{equation}\label{extliecomp}
      A\simeq H\(C^\bdot(L),d\)
  \end{equation}
  with the cohomology algebra of DG algebra $(C^\bdot,\dc^*)$, which
  has $C^m=\L^m L^*=S^m\ovl L^*$\,,
  the differential $\dc^*:C^{m}\rTo C^{m+1}$ is acting as
  \begin{equation}\label{dsd}
     \dc^*\psi(\ovl e_1\ovl e_2\,\cdots\,\ovl e_m)=\!\!\!
     \sum_{1\le i<j\le m}\!\!\!\pm\,\psi(\ovl{[e_i,e_j]}\,
     \ovl e_1\,\cdots\,\without{\ovl e}i\,\cdots\,
                       \without{\ovl e}j\,\cdots\,\ovl e_m)\;,
  \end{equation}
  and the multiplication in $C^\bdot$ is given by the shuffle product\footnote{the summation in
  \eqref{sp} goes over all $I=\{i_1<i_2<\,\cdots\,<i_k\}$,
  $J=\{j_1<j_2<\,\cdots\,<j_m\}$ such that $I\sqcup
  J=\{1,\,2,\,\dots\,,\,(m+k)\}$}
  \begin{equation}\label{sp}
     [\f\comp\psi]\,(\ovl e_1\ovl e_2\,\cdots\,\ovl e_{k+m})=
     \sum\pm
     \f(\ovl e_{i_1}\ovl e_{i_2}\,\cdots\,\ovl e_{i_k})\cdot
     \psi(\ovl e_{j_1}\ovl e_{j_2}\,\cdots\,\ovl e_{j_m})\;.
  \end{equation}
  The isomorphism \eqref{extliecomp}
  takes the internal graded component $A_i$
  to the $i$-th internal degree component of the $i$-th cohomology
  space. It comes from the quasiisomorphism \eqref{A-quism} and
  quasiisomorphic embedding \eqref{s-alt} of Chevalley's resolution
  \eqref{ChC-display} into bar resolution \eqref{bar}.
\end{proposition}

\subsection{Differential perturbation lemma}
  The proof of the main results of the next sect.\;\ref{HcompT} will be
  based on the lemma influenced by A.~Losev's talks on `enhanced
  spectral sequences', which is a slight variation on the simplest,
  degree one, case of the homotopy structure transferring in
  Kadeashvili's type of thing\footnote{comp. with \cite{Bw}; see also
  \cite{Ma} for similar explicit $A_\infty$-formulas most closed to our
  framework}.

  Let $(E\,,\;d:E\rTo E\,,\;d^2=0)$ be an arbitrary differential
  $S$-module over an arbitrary ring $S$. Assume we are given
  with the diagram of $S$-modules and $S$-linear homomorphisms
  \begin{equation}\label{EHdia}
     E\pile{\rTo^{\quad\l\quad}\\\lTo_{\quad\r\quad}}H
  \end{equation}
  together with $S$-linear homotopy $E\rTo^\k E$ that satisfy the
  following properties:
  \begin{gather}
       \l\r=1\;,\qquad\r\l=1+d\k+\k d\;,\label{l-r-rels}\\
       \k^2=d^2=\l d=d\r=\l\k=\k\r=0\label{d-k-rels}
  \end{gather}
 (this means that $H$ can be included into $E$ as the retract of $d$
 capturing all its homology). We intend to show that under this
 assumptions and appropriate `convergency condition' any perturbation
 $D=d+\d$ (even not necessary commuting with $d$) induces some non
 trivial differential $\PD$ on $H$ and a perturbation $(\l',\r',\k')$,
 of the diagram \eqref{EHdia} and the homotopy $\k$, such that $\l'$,
 $\r'$ remain to be the inverse homotopy equivalences between the
 complexes $(E,D)$ and $(H,\PD)$. The convergency condition in question
 is the existence of $\k',\,\e_\l,\,\e_\r\in\End_S(E)$ defined by the
 following series
 \begin{equation}\label{ka-pert}
 \begin{gathered}
    \k'\bydef\k+\k\d\k+\k\d\k\d\k+\k\d\k\d\k\d\k\,+\dots
        =\k(1+\e_\l)=(1+\e_\r)\k\;,\\
     \e_\l\bydef\sum_{m\ge1}(\d\k)^m=\d\k'\;,\quad
     \e_\r\bydef\sum_{m\ge1}(\k\d)^m=\k'\d\;.
  \end{gathered}
 \end{equation}
 Note that these operators are well defined, for example, if $\d\k\in\End_S(E)$ is locally nilpotent\footnote{\ie $\fain eE$ $\exists\;m=m(e)\in\NN$ : $(\d\k)^me=0$}.

\begin{lemma}\label{dpl}
  Under the assumptions \eqref{EHdia}--\eqref{d-k-rels}, let
  $$D=d+\d\,:\;E\rTo E\;,\quad D^2=0\,,
  $$
  be another $S$-linear differential on $E$. If the operators \eqref{ka-pert} are well defined then the perturbed operators\quad $\l'\bydef\l(1+\e_\l)$\;,\quad $\r'\bydef(1+\e_\r)\r$\quad  satisfy the conditions
  \begin{enumerate}
     \item
     $\l'\r'=\Id_H$\;,\quad
     $\r'\l'=\Id_E+D\k'+\k'D$\,;
     \item
     $\PD=\l'D\r=\l'\d\r=\r'D\l=\r'\d\l$ is actually the
     same operator on $H$\,;
     \item
     $\PD^2=0$\,, \ie $\PD$ provides $H$ with the differential;
     \item
     $\PD\l'=\l'D$ and $\r'\PD=D\r'$\,, \ie
     $(E,D)\pile{\rTo^{\quad\l'\quad}\\\lTo_{\quad\r'\quad}}(H,\PD)$
     are morphisms of complexes providing the
     inverse to each other homotopy equivalences.
  \end{enumerate}
\end{lemma}

\begin{proof}
  It follows from \eqref{d-k-rels} that
  $\e_\l\r=\l\e_\r=\e_\l\e_\r=0$\,. This implies
  $$\l'\r'=\l(1+\e_\l)(1+\e_\r)\r=\l\r=\Id_E\;,
  $$
  which is the first relation in (1). The conditions $d^2=0$
  and $D^2=(d+\d)^2=0$ imply that $\d^2=-d\d-\d d$\,. Using this
  relation and \eqref{ka-pert}, we get
  \begin{equation}\label{ereldec}
  \begin{aligned}
     \e_\r\e_\l&=\k(1+\e_\l)\d^2(1+\e_\r)\k=\\
              &=-\k(1+\e_\l)d\d(1+\e_\r)\k
                -\k(1+\e_\l)\d d(1+\e_\r)\k=\\
              &=-\k'd\e_\l-\e_\r d\k'\;.
  \end{aligned}
  \end{equation}
  Now, to compute $\r'\l'$, we substitute $\r\l=1+d\k+\k d$ and
  write the result as a sum of three terms
  \begin{multline}
     \r'\l'=(1+\e_\r)\r\l(1+\e_\l)=(1+\e_\r)(1+d\k+\k d)(1+\e_\l)=\\
     =(1+\e_\r)(1+\e_\l)+
      (1+\e_\r)d\k(1+\e_\l)+
      (1+\e_\r)\k d(1+\e_\l)\;,
  \end{multline}
  then expand these summands using \eqref{ereldec} and
  \eqref{ka-pert}\;:
  $$\begin{gathered}
      (1+\e_\r)(1+\e_\l)=1+\e_\r+\e_\l+\e_\r\e_\l=
       1+\d\k'+\k'\d-\e_\r d\k'-\k'd\e_\l\;,\\
    (1+\e_\r)d\k(1+\e_\l)=(1+\e_\r)d\k'=d\k'+\e_\r d\k'\;,\\
   (1+\e_\r)\k d(1+\e_\l)=\k'd(1+\e_\l)=\k'd+\k'd\e_\l\;.
    \end{gathered}
  $$
  Adding up the right sides, we get the homotopy relation required in (1)
  $$\r'\l'=1+\d\k'+\k'\d+d\k'+\k'd=1+D\k'+\k'D\;.
  $$
  Since we have $\l D=\l\d$, $D\r=\d\r$,
  $\e_\l\d=\d\e_\r$\,, (2) follows
  $$\l'D\r=\l'\d\r=\l(1+\e_\l)\d\r=\l\d(1+\e_\r)\r=\l\d\r'=\l D\r'\;.
  $$
  The commutation relations (4) also follow from  $\l D=\l\d$\,,
  $D\r=\d\r$ using (1)
  $$\begin{aligned}
      \PD\l'=
      \l D\r'\l'=\l D(1+D\k'+\k'D)=\l(1+\d\k')D=\l'D\;,\\
      \r'\PD=\r'\l'D\r=(1+D\k'+\k'D)D\r=D(1+\k'\d)\r=D\r'\;.
    \end{aligned}
  $$
  Finally, $\PD\PD=\l D\r'\PD=\l D^2\r'=0$ gives (3).
\end{proof}

\subsection{\mathversion{bold}Chevalley's complex as $S(V^*)$-module}\label{HcompT}
   Consider graded Lie ideal
   $$\Lt=\OP_{m\ge2}L_m\subset L
   $$
   and denote
   by $(C^\bdot(\Lt),\dc^{\ge 2})$ its Chevalley complex. The
   s-exterior algebra of $L$ splits as the graded algebra into the
   tensor product
   \begin{equation}\label{CD-dec}
       C^\bdot(L)=\L(\Lt^*\oplus L^*_1)=\L(\Lt^*)\OX_\CC\L(L_1^*)=
       C^\bdot(\Lt)\OX_\CC S
   \end{equation}
   where\footnote{recall that $L_1=V$ has internal degree 1, thus,
   its s-exterior algebra is nothing but the ordinary symmetric
   algebra} $S=S(V^*)=\L(L_1^*)$ is the projective coordinate
   algebra of $\PP(V)$. The both sides of \eqref{CD-dec} carry the
   natural structure of right $S$-modules\footnote{since $S$ is a
   commutative algebra, it does not matter from what side does it
   act from, but we use the right action to outline that it
   commutes with the left $B$-action on the cobar complex} coming
   from the algebra inclusion
   $$S(V^*)=\L(L_1^*)\rInto\L(L^*)
   $$
   and the isomorphism \eqref{CD-dec} is clearly $S$-linear.
   Moreover, the Chevalley differential \eqref{dsd}, acting on the
   left side, is also $S$-linear because of $\dc(L_1)=0$. The right
   side of \eqref{CD-dec} carries the intrinsic $S$-linear
   differential
   \begin{equation}
      d=\dc^{\ge2}\ox 1
   \end{equation}
   induced by the Chevalley differential for the Lie s-algebra
   $\Lt$. We transfer it to the left side preserving the notation
   $d$ for it. Then on the left side we have
   \begin{equation}\label{CH-d-de}
       \dc=d+\d
   \end{equation}
   where $d$ acts on the subalgebra $\L(\Lt^*)$ by the same formula
   \eqref{dsd} and annihilates all the monomials containing
   $L_1^*$-factors, and $\d$ is the difference, which is automatically
   $S$-linear as well. We are in a position to apply the
   differential perturbation lemma of \ref{dpl}.

\begin{theorem}\label{th3-1}
   Let $A$ be commutative Koszul quadratic algebra,
   $$B=A^!=U(L)
   $$
   be its dual, treated as the universal enveloping algebra of
   the graded Lie s-algebra $L$ (see \ref{QDA}). Then  for
   each $p\ge1$ and any $q$ there exists an isomorphism
   \begin{equation}\label{HcompT-eq}
       R_{p,q}(A)\simeq H^{q-p}(\Lt,\CC)_q
   \end{equation}
   between $q$-th internal degree components of $p$-th syzygy space
   of $A$ {\upshape(}see \eqref{Rpdef}{\upshape)} and $(q-p)$-th
   cohomology space of $\Lt$.
\end{theorem}

\begin{proof}
   Let us split $C^\bdot(\Lt)$ as the vector space over $\CC$
   into a direct sum of bigraded subspaces
   \begin{equation}\label{SD-C}
      C^\bdot(\Lt)=H\oplus I\oplus P
   \end{equation}
   where $I=\im\dc^{\ge2}$ and $H\oplus I=\ker\dc^{\ge2}$. Thus
   $H\simeq H^\bdot(\Lt)$ and $d$ takes $P$ isomorphically onto $I$
   and annihilates $H\oplus I$. Write $\k$ for the operator
   \begin{equation}\label{kapL-def}
      H\oplus I\oplus P\rTo^\k H\oplus I\oplus P
   \end{equation}
   that annihilates $H\oplus A$ and acts on $I$ as
   $-d^{-1}:I\rTo^\sim P$. We write
   \begin{equation}\label{HIP-dia-L}
     C^\bdot(\Lt)\pile{\rTo^{\quad\l\quad}\\\lTo_{\quad\r\quad}}H
   \end{equation}
   for the embedding and the projection associated with the direct
   sum decomposition \eqref{SD-C}. So, our $\k$, $\l$, $\r$
   satisfy the relations \eqref{l-r-rels}, \eqref{d-k-rels}.
   Tensoring \eqref{HIP-dia-L} by $S$ and combining it with the
   $S$-module isomorphism \eqref{CD-dec}, we get the diagram of
   $S$-modules
   \begin{equation}\label{HIP-dia-LS}
     (C^\bdot(L),\dc(L))\pile{\rTo^\l\\\lTo_\r}
     H^\bdot(\Lt)\OX_\CC S\;,
   \end{equation}
   and the $S$-linear map $C^\bdot(L)\rTo^\k C^\bdot(L)$ that
   satisfy the relations \eqref{l-r-rels}, \eqref{d-k-rels}
   as well. Further, the composition $\d\k$, where $\d$ comes
   from the decomposition \eqref{dsd}\,, is locally a nilpotent
   operator, because it is clear from \eqref{dsd} that $\d d^{-1}$
   preserves the homological degree and strictly decreases the
   difference between the total internal degree and degree induced
   by the homological grading coming from $C^\bdot(\Lt)$. Thus, the
   differential perturbation lemma from \ref{dpl} provides
   $H^\bdot(\Lt)$ with the differential
   \begin{equation}\label{perturebedLD}
       \PD=\l\comp\(\sum_{m\ge0}(\d\k)^m\)\comp\r
   \end{equation}
   such that $(H^\bdot(\Lt)\ox S,\PD)$ becomes a complex of free graded
   $S$-modules {\it homotopy equivalent\/} to the Chevalley complex
   $(C^\bdot(L),\dc(L))$. Since the latter is quasiisomorphic to $A$
   as a DG $S$-module, we can compute
   $$R_p=\Tor_p^S(A,\CC)
   $$
   as $p$-th cohomology of the complex obtained by
   tensoring $(H^\bdot(\Lt)\ox S,\PD)$ by $\CC$ over $S$. Because the
   differential \eqref{perturebedLD} strictly increases the internal
   $S$-module degree, it will be annihilated by this tensoring and
   we get the required isomorphism \eqref{HcompT-eq}.
\end{proof}

\subsubsection{Coincidence of two algebra structures on the syzygies.}
  In \S2 we have equipped the space of syzygies by an
  algebra structure induced by the Koszul DGA resolution
  \eqref{cAalgebra} for $A$. On the other side, the Lie algebra
  cohomology $H^\bdot(\Lt)$ also has an algebra structure induced
  by the Chevalley's DGA resolution. So, the both sides of
  \eqref{HcompT-eq} come with the intrinsic algebra structures. In
  fact this two structures do coincide.

\begin{theorem}
   The isomorphism $R\simeq H^\bdot(\Lt)$ constructed in
   the previous theorem is an isomorphism of algebras.
\end{theorem}

\begin{proof}
  We have to compare Chevalley's DG algebra $C^\bdot(\Lt)$ with the
  DG algebra
  \begin{equation}\label{A-koz2}
     \cA=A\OX_S K_S\;,
  \end{equation}
  where $K_S=S^\bdot(V^*)\OX_\CC\L^\bdot(V^*)$ is the Koszul
  complex \eqref{KRforC}, which is the Koszul resolution for $\CC$
  as the left module over $S=S(V^*)$. To this aim consider
  \begin{equation}\label{bigDGA}
      E=C^\bdot(L)\OX_S K_S=C^\bdot(L)\OX_\CC\L(V^*)
  \end{equation}
  equipped with the differential $D=\dc\OX_S1+1\OX_S\dk$, where
  $\dk$ is the Koszul differential \eqref{KD}. This DG
  algebra has the compatible structure of right DG module over
  DG algebra $K_S$. The $S$-module isomorphism \eqref{CD-dec}
  extends obviously to the isomorphism of right $K_S$-modules
  \begin{equation}\label{K-isom}
       C^\bdot(L)\OX_SK_S=E\;\simeq\;E'=C^\bdot(\Lt)\OX_\CC K_S\;.
  \end{equation}
  The right side has the intrinsic $S$-linear differential
  $d=\dc^{\ge2}\ox1+1\ox\dk$, which can be transferred to $E$.
  Thus, we get two decompositions for $D$ into a sum
  \begin{equation}\label{twoDdec}
     D=\dc\OX_S1+1\OX_S\dk=d+\d
  \end{equation}
  We are going to compute the algebra structure on the cohomology
  space $H(E,D)$ using the spectral sequences associated with these two
  decompositions.

  The first decomposition $D=\dc\ox1+1\ox\dk$ has commuting
  summands, \ie represents $E$ as a double complex. Those spectral
  sequence that firstly computes the Koszul cohomology degenerates
  in $E_2$-term. Its $E_1$-term is concentrated at the zero row and
  coincids with the Koszul resolution \eqref{A-koz2}. Thus,
  $H(E,D)\simeq R$ and has the multiplication induced from
  \eqref{A-koz2}.

  The second decomposition $D=d+\d$ corresponds to the natural
  filtration on $E$ coming from
  the Serre--Hochschild filtration of the pair $(L,\Lt)$ on
  $C^\bdot(L)$, where
  $q$-th filtered component is dual to the $\CC$-linear span of all
  monomials \eqref{Lbmons} that contain $\le q$ generators
  $e_i\in\ovl L_1\subset\ovl L$. It is clear from \eqref{dsd},
  \eqref{sp} that this filtration is compatible with $D$ and the
   product in $E$. Thus, we conclude that the algebra structure
  on $H(E,D)$ coincides with the one induced from the DG algebra
  structure on the $E_1$-term of the spectral sequence for this filtration.

  Obviously this $E_1$-term is nothing but the right side of \eqref{K-isom} with its natural algebra structure and differential $d=\dc^{\ge2}\ox1+1\ox\dk$. To compute its homology we can use the fact that it, in its own turn, is the sum of two commuting differentials. Applying the same arguments as above, we conclude that the $E_2$-term is isomorphic to  $H^\bdot(\Lt)$  and has the algebra structure induced by the Chevalley DGA resolution for $\Lt$. Now Theorem \ref{th3-1} implies that this spectral sequence also degenerates at  $E_2$-term. We conclude that  there is an algebra isomorphism $H(E,D)\simeq H^\bdot(\Lt)$. Thus $ R  \simeq H(E,D)\simeq H^\bdot(\Lt)$ and we can say that the multiplicative structures on $R$ is induced from $C^\bdot(\Lt)$.
\end{proof}

\section{Some geometric examples}

\subsection{Notation and preliminaries}\label{GLnot}
  In this section we use Theorem \ref{th3-1} to describe the syzygies
  \begin{equation}\label{rpq-hq-p-q}
     R_{p,q}(A)\simeq H^{q-p}(\Lt,\CC)_q
  \end{equation}
  of projective coordinate algebras of certain HW-orbits in $\PP(V)$ by computing the cohomologies staying in the right hand side of \eqref{rpq-hq-p-q}. Recall (see \ref{SHWex-2}) that the bigraded components  \eqref{rpq-hq-p-q} the syzygies of an HW-orbit are equipped with the natural action of $\GL(V)$. Our computations will use irreducible decompositions of \eqref{rpq-hq-p-q} \wrt this action.

\subsubsection{Young diagram notations.}
  We depict a partition\footnote{our notations for the Young tableaux and associated symmetric functions agree
  with \cite{Fu2}, \cite{MD} where the reader can find all the
  formulas we will use below to express the symmetric polynomials
  through each other} $\l=[\FAM\l,k]$
  (where $\l_1\ge\l_2\ge\dots\ge\l_k$) by Young diagram with
  $k$ rows of lengths $\FAM\l,k$ like
  \begin{equation}\label{YDsample}
     \hbox{\tiny$\yng(5,3,3,1)$}\qquad(\text{for}\quad\l=[5,3,3,1])
  \end{equation}
  and write $\l'=[\FAM\l',m]$ for the transposed diagram (say, for \eqref{YDsample} we have $\l'=[4,3,3,1,1]$). The total number of cells $|\l|=\sum\l_i$ will be called {\it a weight\/} of the diagram. The shorted notation $[\r_1^{s_1},\,\r_2^{s_2},\dots,\,\r_n^{s_n}]$ means the Young diagram that has $s_i$ rows of length $\r_i$ (in  \eqref{YDsample} $\l=[5,3^2,1]$).

  Also we will use the Frobenius notation and write
  \begin{equation}\label{Fnot}
     (\FAM\a,p\,|\,\FAM\b,p)
  \end{equation}
  for the Young diagram $\l$ whose main diagonal consists of $p$
  cells and $\l_i=\a_i+i$\,, $\l'_i=\b_i+i$ for each $i=1,\,2,\,\dots\,,\,p$
  (in  \eqref{YDsample} $\l=(4,1,0\,|\,3,1,0)$). Note that in this
  notation $\a_1>\a_2>\cdots>\a_p\ge0$ and
  $\b_1>\b_2>\cdots>\b_p\ge 0$.

\subsubsection{Irreducible $\GL$-modules.}\label{GLNchars}
  With any Young diagram
  $\l=[\FAM\l,k]$
  is associated an irreducible
  $\GL_k(\CC)$-module $\pi_\l$ of the highest weight
  $$\l_1\e_1+\l_2\e_2+\dots+\l_k\e_k
  $$
  which is simultaneously the irreducible $\SL_k(\CC)\subset\GL_k(\CC)$\;-\;module of highest weight
  $$(\l_1-\l_2)\,\a_1+(\l_2-\l_3)\,\a_2+\dots+(\l_{k-1}-\l_k)\,\a_{k-1}
  $$
  (where $\a_i=\e_i-\e_{i+1}$ are the simple roots as in
  \ref{g25}). Recall that $\pi_\l$ can be constructed by factorizing the space
  \begin{equation}\label{eq:La-space}
     \L^\l\bydef\L^{\l'_1}V\ox\L^{\l'_2}V\ox\;\cdots\;\ox\L^{\l'_m}V
  \end{equation}
  (where $V=\CC^k$ is the tautological $\GL_k$-module and the exterior powers are the lengths of columns of $\l$) through {\it the column exchange relations\/}\footnote{see \cite{Fu2}}. Such a relation can be written for any choice of the following data:
  \begin{itemize}
    \item
    a filling $T$ of the Young diagram $\l$ by $n=|\l|$ vectors $\FAM v,n\in V$
    \item
    a number $i$ of a column in $\l$
    \item
    a collection $I$ of cells in the next $(i+1)$-th column
  \end{itemize}
  Taking an exterior product of the vectors $v_\nu$ along each column of $\l$ and tensoring these products together, we get an element $v^T\in\L^\l$; then the column exchange relation, which corresponds to $T,i,I$ says that in $\pi_\l$
  \begin{equation}\label{eq:exchangerel}
     v^T=\sum_\s v^{\s T}\;,
  \end{equation}
  where $\s$ runs through the permutations of $v_\nu$'s providing order preserving exchanges between $I$-cells and all collections of $\#I$ cells in the previous $i$-th column\footnote{if $\#I=m$ there are totally $\binom{\l'_i}{m}$ such permutations}. For example, an exchange relation corresponding to a filled diagram
  \hbox{\tiny$\young(13,2)$} and the only possible choice $i=1$, $I=\{2\}$ says that  $\forall\;v_1,v_2,v_3$ in $\pi_{(2,1)}$ we have
  $$(v_1\w v_2)\ox v_3=(v_3\w v_2)\ox v_1+(v_1\w v_3)\ox v_2\;,
  $$
  which is nothing but the Jacobi relation
  $[[v_1,v_2],v_3]+[[v_2,v_3],v_1]+[[v_3,v_1],v_2]=0$\,.

  If we write $\FAM e,k$ for the standard basis in $V=\CC^k$, then the standard basis for $\pi_\l$ is formed by classes of elements $e^T\in\L^\l$ obtained from fillings $T$ of $\l$ by vectors $e_i$ such that the indexes of $e_i$ weakly increase across each row and strictly increase down each column, \ie form {\it a Young tableau\/} $T$ of the shape $\l$ on the alphabet $[1\,.\,.\,k]$. With any such a tableau $T$ one can associate a monomial $x^T=x_1^{m_1}x_2^{m_2}\dots x_k^{m_k}$\,, where $m_i$ is number of occurrences of $i$ in $T$. Then the character of the irreducible $\GL_k$-module $\pi_\l$ is the Schur polynomial
  $$s_\l\VEC x,k=\sum_Tx^T\;,
  $$
  where the sum is running over all Young tableaux.
  For example, a filling
  $$T=\hbox{\tiny$\young(11125,334,445,5)$}
  $$
  of the diagram \eqref{YDsample} is a valid tableau for $\GL_5$ and contributes monomial
  $x_1^3x_2x_3^2x_4^3x_5^2$ into $s_{[5,3,3,1]}$\,: this monomial
  computes the eigenvalue of the standard basic vector coming from
  $$e^T=\(e_1\w e_3\w e_4\w e_5\)\ox
    \(e_1\w e_3\w e_4\)\ox
    \(e_1\w e_4\w e_5\)\ox
      e_2\ox
      e_5\,.
  $$

\subsubsection{Euler's $\GL$-characters.}\label{Echars}
  Since the category of $\GL$-modules is semisimple, each term of
  any $\GL$-equi\-va\-ri\-ant complex $K^\bdot$ splits as
  $$K^\nu=\OP_{\l\in\L_\nu}\pi_\l\,,
  $$
  where $\L_\nu$ is the set of highest weights of all irreducible
  representations appearing in $K^\nu$ (counted with
  multiplicities). We call the alternated sum
  \begin{equation}\label{GLchi}
     \chi_{K^\bdot}\bydef\sum_\nu(-1)^\nu\sum_{\l\in\L_\nu}s_{\l}
  \end{equation}
  {\it the Euler $\GL$-character\/} of $K^\bdot$. Clearly, $\chi_{K^\bdot}=\chi_{H(K^\bdot)}$, where $H(K^\bdot)$ is considered as a complex with zero differentials. Similarly, one can define the Euler $\GL$-characters for graded $\GL$-equi\-va\-ri\-ant commutative and Lie s-algebras. A powerful tool for comparing these characters is provided by the following version of the Koszul duality.

\subsubsection{Quillen's duality.}
In discussion before Proposition \ref{CH-summary} we associated the Chevalley complex $C^\bdot(\ga)$ with any Lie s-algebra $\ga$. It can be treated as commutative s-algebra $\L^\bdot(\ga^*[-1])$ equipped with the differential whose action on generators is dual to the bracket $\L^2\ga\rTo\ga$. This construction provides a functor from the category of Lie s-algebras to one of commutative DG-algebras.

On the other side, with any commutative s-algebra $A$ one can associate in a similar way the Harrison complex $\cH^{\bdot}(A)$, which is a free Lie s-algebra $\Lie^{\bdot}(A^{*}[-1])$ equipped with the differential whose action on generators is dual to the multiplication $S^2A\rTo A$. Thus, we get a functor acting in the opposite direction.

Applying to these functors the general result proven in \cite[th.\,4.2.5]{GK} for any pair of Koszul dual operads, we get

\begin{lemma}
  Two functors described above are homotopy inverse to each other, \ie
  for any pro-nilpotent Lie s-algebra $\ga$ and any s-commutative
  algebra $A$ there exist natural quasi-isomorphisms
  $C^\bdot(\cH^\bdot(A))\rTo^\sim A$\,,\quad
  $\cH^\bdot(C^\bdot(\ga))\rTo^\sim\ga$\,.\EPbox
\end{lemma}

\begin{corollary}\label{FLch}
   Let $\ga$ be a pro-nil\-po\-tent graded Lie s-al\-geb\-ra
   equipped with $\GL_k$-ac\-tion preserving the graded Lie s-al\-geb\-ra
   structure. Then the Euler $\GL_k$-cha\-rac\-ters of $\ga$ and
   $H=H^\bdot(\ga,\CC)$ are related by
   \begin{align}
     \chi_H&=
     \sum_{n\ge0}(-1)^{n}e_n\comp\chi_\ga\;,\label{eq:LieA2Hch}\\
     \chi_\ga&=
     \Bigl(-\sum_{m\ge1}\mu(m)\ln\bigl(1-(-1)^mp_m\bigr)\Bigr)\comp
     \chi_H\;,\label{eq:LieH2Ach}
   \end{align}
  where $\mu$ is the M\"obius function, $e_n$ are the elementary
  symmetric polynomials, $p_m$ are Newton's sums of powers, and
  $\comp$ means the plethysm of symmetric functions.
\end{corollary}

\begin{proof}
It is well known that $\GL$-character of $\L^{n}(V)$ equals $e_n$ (see \cite{MD}) and $\GL$-character of $\Lie^{n}(V)$ equals $\frac{1}{n}\sum_{d|n}\mu(d)p_{d}^{\frac{n}{d}}$ (see \cite{Klyachko}). The signs in \eqref{eq:LieA2Hch} and \eqref{eq:LieH2Ach} come from the grading shift.
\end{proof}

\subsubsection{Free Lie algebras.}\label{fLee}
  If $\ga=\Lie(W)$ be a free Lie (graded, s-) algebra generated by
  a vector space $W$, then its universal enveloping algebra
  $U(\ga)$ is a free associative algebra generated by $W$ and the
  trivial $U(\ga)$-module $\CC$ admits a short free resolution
  \begin{equation}\label{resol_k}
     0\rTo U(\ga)\ox W\rTo U(\ga)\rTo\CC\rTo0\;.
  \end{equation}
  For any subalgebra $\gb\subset\ga$ an isomorphism of
  $\gb$-modules $U(\ga)\simeq U(\gb)\ox S(\ga/\gb)$ shows that
  \eqref{resol_k} is a free resolution for $\CC$ in the category of
  $U(\gb)$-modules as well. It is well known that pro-nilpotent Lie
  algebra $\gb$ is free iff $H^i(\gb,\CC)=0$ for $i>1$. In
  particular, applying $\hom_{U(\gb)}(\,\ast\,,\CC)$ to the short
  resolution \eqref{resol_k}, we get

\begin{lemma}\label{SUBfLee}
  A subalgebra of any free Lie (super)algebra is free.\EPbox
\end{lemma}

\subsection{The most singular case}
  To begin with, consider the commutative quadratic algebra with the maximal possible space of quadratic relations, \ie
  $$A=S(V)/(S^2V)\;.
  $$
  Its Koszul dual Lie s-algebra is the free graded Lie algebra $L=\Lie\(V[-1]\)$ generated by the vector space $V$ situated in degree 1. The existence of resolution \eqref{resol_k} (written for $\ga=L$) and criteria from sect.\;\ref{KozCrit} imply that $L$ and $A$ are Koszul. By Lemma \ref{SUBfLee} the Lie subalgebra $\Lt\subset L$ is free as well. This implies that the syzygies \eqref{rpq-hq-p-q} vanish for $(q-p)>1$. Thus, non trivial syzygies are described by

\begin{proposition}\label{syz_free_Lie}
   For $A=S(V)/S^2(V)$ the component $R_{p,(p+1)}$ of
   syzygies \eqref{rpq-hq-p-q} is the irreducible $\GL(V)$-module $\pi_{[2,1^{p-1}]}$. All the other syzygies vanish.
\end{proposition}

\begin{proof}
  In this case the Koszul complex \eqref{KCforA} takes extremely  simple form and can be written as the tensor product
  $K^{\bdot}=(\CC\oplus V)\ox \L^{\bdot}(V[-1])$. Thus, its Euler $\GL(V)$-character has the following expression in terms of the elementary symmetric polynomials $e_k$
  \begin{equation}
     \chi_{K^\bdot}=
     (1+ e_1)\cdot\!\sum_{k\ge0}(-1)^k\,e_k=
     1+\sum_{k\ge1}(-1)^k\,(e_1e_k-e_{k+1})\;.
  \end{equation}
  By the Frobenius character formula\footnote{see, for example,
  \cite[ch.\,1.3,\;formula\;(3.5)]{MD}}, the latter multiplier of degree $(k+1)$
  coincides with the irreducible $\GL(V)$-character
  $$s_{[2,1^{k-1}]}=\det\begin{pmatrix}e_k&e_{k+1}\\e_{0}&e_1\end{pmatrix}
  $$
  of the irreducible $\GL$-module $[2,1^{k-1}]$. Since for each $q>0$
  there is only one non zero component $R_{p,q}$ and it has
  $p=q-1$, we conclude that $[2,1^{k-1}]=R_{k,(k+1)}$.
\end{proof}

\subsection{Syzygies of the Veronese curve}
  The Veronese embedding takes
  $$\PP_1=\PP(U)\rInto^{\quad u\mapsto u^d\quad}\PP(S^nU)=\PP_n\;.
  $$
  In appropriate coordinates on $\PP(S^nU)$ it sends
  $(u_0:u_1)\in\PP_1$ to
  $$(x_0:x_1:\,\dots\,:x_n)=
    (u_0^n:u_0^{n-1}u_1:u_0^{n-2}u_1^2:\,\dots\,:u_1^n)\;.
  $$
  The image is described by the quadratic equations $x_ix_j=x_kx_m$ (for all possible choices of $i$, $j$, $k$, $m$ with $i+j=k+m$), so, we get
  the quadratic algebra
  $$A=\CC[x_0,\,\dots\,,\,x_n]/\!\(x_ix_j-x_kx_m\,|\;i+j=k+m\)\;.
  $$
  In the dual coordinates $x^i$ on $V^*=S^nU^*$ the relations for the
  Koszul dual Lie s-algebra take a form
  \begin{equation}\label{VerLrel}
    \sum_{i+j=k}[x^i,x^j]=0\;,\;\text{ where }k=0,\,1,\dots,2n\;.
  \end{equation}
  Let us order the generators in the following non-standard way
  $$x^0>x^n>x^{n-1}>\dots>x^1\;.
  $$
  This ordering induces the filtration on $L$. Since the leading Lie monomial of $k$-th relation in \eqref{VerLrel}
  is $[x^0,x^k]$ for $k\le n$ and is $[x^n,x^k]$ for $k\ge n$, associated with this filtration graded Lie algebra $L'$ is isomorphic to the direct sum of the abelian Lie algebra generated by $x^0,\,x^n$ and the free Lie algebra $\Lie(W)$, where $W$ is the linear span of remaining $(n-1)$ variables $x^1,\dots,x^{(n-1)}$. Since the abelian part makes no contribution in the Chevalley complex for $\Lt$, we have $C^\bdot(\Lt)=C^\bdot\bigl(\Lie(W)_{\ge2}\bigr)$ and can repeat word for word the computations of the previous section taking $W$ instead of $V$.

  \begin{corollary}\label{cor:VCsyz}
      The only non zero syzygies of the Veronese curve in $\PP_n$ are $R_{p,p+1}=\pi_{[2,1^{p-1}]}$ with\footnote{we consider $R_{p,p+1}$ as $\GL(W)$-modules and apply the hook length formula (see \cite{Fu1}, \cite{MD})} $\dim R_{p,p+1}=p\cdot\binom n{p+1}$. \EPbox
  \end{corollary}

\subsection{\mathversion{bold}
Syzygies of the grassmannian $\Gr(2,N)$}\label{answer_grasm}
  This computation generalizes the examples \ref{g25}, \ref{SHWex-2c}.
  The Pl\"ucker embedding
  \begin{equation}\label{G2Npluemb}
     \Gr(2,N)\rTo^{\sim}\PP(G\cdot v_{\hw})\subset\PP(V_{\o_2})
  \end{equation}
  realizes the grassmannian as the \SHW of $G=\SL(N,\CC)$ with the index \eqref{QHWeq} equal to $N$. Moreover, the natural action of $\GL_N(\CC)$ on $\PP(V_{\o_2})$ preserves the grassmannian and is kept on the syzygies as well. Write $W=V_{\o_1}^*$ for dual to the standard $N$-dimensional representation of $\GL_N$. It follows from Propositions \ref{qe-proof-pre} and \ref{qe-proof} that the projective coordinate algebra of the  Pl\"ucker embedding \eqref{G2Npluemb} has the following $\GL_N$-module decomposition:
  \begin{equation}\label{G2NAintro}
     A=S^\bdot(\L^2 W)/(\L^4W)=\OP_{k=0}^{\infty}\pi_{[k,k]}\;.
  \end{equation}
  Our computation of the syzygies of $A$ will be organized as follows.
  In section \ref{G2Nstep1} we compute the Euler \GLc of the Koszul complex
  \begin{equation}\label{G2NKintro}
     K^\bdot=\L^{\bdot}\((\L^2 W)[-1]\)\ox A\;.
  \end{equation}
  Then, in section \ref{G2Nstep2}, we completely describe the generators of $\Lt$, \ie compute $H^1(\Lt,\CC)$. It turns out that these generators form a graded $\GL$-module\footnote{we are using here the Frobenius notations \eqref{Fnot}}
  $$\OP_{2\le q\le(N-2)}\hspace{-1ex}
    \pi_{((q-2)\,|\,(q+1))}\quad\text{with}\quad
    \pi_{((q-2)\,|\,(q+1))}=R_{q-1,q}=H^1(\Lt,\CC)_{q}\;.
  $$
  Moreover, we will show that $\pi_{\((q_1-1),\dots,(q_p-1)\,|\,(
  q_1+2),\dots,(q_p+2)\)}$ appears with multiplicity one in the syzygy space
  $R_{p,(q_1+\dots+q_p)}$.

  This allows to guess the shape of answer and forces to introduce a bigraded skew commutative $\GL_N$-equivariant s-algebra
  \begin{equation}\label{eq:fAdef}
  \begin{gathered}
    \fA=\OP_{p,q}\fA_{p,q}\;,\quad\text{where}\\
    \fA_{p,q}=\hspace{-1ex}\OP_{\substack{
    (N-2)\ge i_1>\dots>i_p\ge2\\i_1+\dots+i_p=q}}\hspace{-1ex}
    \pi_{((i_1-2),\dots,(i_p-2)\,|\,(i_1+1),\dots,(i_p+1))}\;.
  \end{gathered}
  \end{equation}
  We define $\fA$ as a skew commutative s-algebra generated by the graded vector space
  \begin{equation}\label{G2NA1}
     \fA_1=\OP_q\fA_{1,q}\,,\text{ where }
      \fA_{1,q}=
    \begin{cases}
       \pi_{(q-2|q+1)}\,,&\text{for $2\le q\le(N-2)$}\,,\\
       0\,,&\text{otherwise}
    \end{cases}
  \end{equation}
  By the definition, the multiplication map $\fA_{1,q_1}\ox\dots\ox\fA_{1,q_p}\rTo\fA_{p,(q_1+\dots+q_p)}$ is given by the projection onto irreducible component
  \begin{equation}\label{eq:fAproj}
     \pi_{((q_1-2),\dots,(q_p-2)\,|\,(q_1+1),\dots,(q_p+1))}\subset
     \pi_{(q_1-2|q_1+1)}\ox\dots\ox\pi_{(q_p-2|q_p+1)}\;,
  \end{equation}
  if $q_1>\dots>q_p$, and vanishes when some of $q_i$'s coincide. Since the component \eqref{eq:fAproj} has the multiplicity one, the corresponding projection is unique up to proportionality and we actually get well defined associative algebra $\fA$ with the components \eqref{eq:fAdef}. This algebra gives a particular example of $\GL$-equivariant {\it hook algebra\/}.

  General properties of hook algebras will be discussed systematically in the next \S5. In particular, in section \ref{G2Nstep3} we show that $\fA$ is quadratic and Koszul.

  Using this result, we show in section \ref{G2Nstep4} below that a
  Lie s-algebra $\fL$, Koszul dual to $\fA$, is isomorphic to $\Lt$.
  This implies the coincidence $\fA=H(\fL,\CC)=H(\Lt,\CC)=R$\,.

\begin{theorem}\label{syzygies}
     The syzygies of the Grassmannian $\Gr(2,N)$ form a bigraded
     skew commutative Frobenius quadratic Koszul algebra isomorphic
     to the algebra $\fA$ defined in
     \eqref{eq:fAdef}--\eqref{eq:fAproj}:
     $$R_{p,q}=\fA_{q-p,q}\,.
     $$
\end{theorem}

\subsection{Koszul complex.}\label{G2Nstep1}
  Euler's \GLc of $K^\bdot=\L^{\bdot}\((\L^2 W)[-1]\)\ox A$ is
  \begin{equation}\label{KGrECpre}
     \chi_{K^\bdot}=
     \chi_{\L^\bdot(\L^2W[-1])}\cdot\Bigl(\sum_js_{[j,j]}\Bigr)=
     \Bigl(\sum_k(-1)^ke_k\comp e_2\Bigr)\;\cdot\;
     \Bigl(\sum_jh_j^2-h_{j+1}h_{j-1}\Bigr)\;,
  \end{equation}
  where $h_k$ are the complete symmetric functions, $\comp$ denotes the plethysm of symmetric functions, and $$s_{[j,j]}=\det\begin{pmatrix}h_j&h_{j+1}\\h_{j-1}&h_j\end{pmatrix}
  $$
  by the Frobenius character formula\footnote{see \cite{Fu2} or \cite[ch.\,1.3,\;formula\;(3.4)]{MD}}. The right hand side of \eqref{KGrECpre} has the following expansion in terms of Schur polynomials:

\begin{lemma}\label{gr_char}
  $\displaystyle
   \chi_{K}=\hspace{-6ex}\sum_{\substack
   {p\ge 0\\(N-3)\ge i_1>\dots>i_p>0}}\hspace{-6ex}
    (-1)^{i_1+\dots+i_p}
    s_{\((i_1-1),\dots,(i_p-1)\,|\,( i_1+2),\dots,(i_p+2)\)}$\,.
\end{lemma}

\begin{proof}
  We start from the generalized
  Littlewood formula established in \cite[Th.\,4.4,\;eq.\,(2)]{IW}:
  \begin{multline}\label{IW-form}
    \sum_{\substack{p\ge0\\n\ge i_1>\dots>i_p\ge0}}\hspace{-4ex}
    (-1)^{p+i_1+\dots+i_p}s_{\((i_1+r),\dots,(i_p+r)\,|\,i_1,\dots,i_p\)}
    =\\[-3ex]
    \prod_{i\le j}^{n}(1-x_i x_j)\;\cdot\;
    \prod_{i=1}^{n}x_i^{\frac{(r-1)}{2}}\cdot
    \frac{\det
    \(x_i^{j+\frac{r-1}{2}}-x_i^{-j-\frac{r-1}{2}}\)_{1\le i,j\le n}}
    {\D_{C(n)}}\;,
  \end{multline}
  where the Weyl determinant $\D_{C(n)}$ in the denominator
  can be expressed as
  \begin{equation}\label{Wdet}
  \begin{gathered}
    \D_{C(n)}\VEC x,n=
    \det\(x_i^{j}- x_i^{-j}\)_{1\le i,j\le n}=\\
    =\prod_{i=1}^n x_i^{-n}\;\cdot\;
     \prod_{i=1}^n(x_i^2-1)\;\cdot\;
     \prod_{i<j}^n(x_i-x_j)(1-x_i x_j)\;.
  \end{gathered}
  \end{equation}
  We take $r=3$ and note that
  $\det\(x_i^{j+1}- x_i^{-j-1}\)_{1\le i,j\le n}$
  staying in the numerator of \eqref{IW-form} can be written as
  \begin{equation}\label{r1det}
    \D_{C(n)}\;\cdot\;
    \prod_{i=1}^n x_i^{-1}\;\cdot\;
    \Bigl(\sum_{j} e_j^2-e_j e_{j+2}\Bigr)\;.
  \end{equation}
  Indeed, consider the $(n+1)$-th order Weyl determinant
  $\D_{C(n+1)}(x_1,\dots,x_n,q)$ as a Lourant polynomial in $q$
  and compute its coefficient at $q^1$ in two ways: by the
  straightforward expansion of the determinant in the middle of
  \eqref{Wdet} and using triple product expansion from the last
  term of \eqref{Wdet}. In the first case we get
  $$(-1)^n\det\(x_i^{j}- x_i^{-j}\)_{\substack
    {1\le i\le n\\2\le j\le n+1}}=
    (-1)^n\det\(x_i^{j+1}- x_i^{-j-1}\)_{1\le i,j\le n}\;
  $$
  In the second case we have
  \begin{gather*}
    \prod_{i=1}^nx_i^{-n-1}
    \prod_{i=1}^n(x_i^2-1)
    \prod_{i<j}^n(x_i-x_j)
    (1-x_i x_j)q^{-n-1}(q^2-1)\prod_{i=1}^{n}(x_i-q)(1-x_i q)=\\
    =\D_{C(n)}\;\cdot\;
    \prod_{i=1}^nx_i^{-1} q^{-n-1}(q^2-1)\;\cdot\;
    \Bigl(\sum_{i=0}^n(-1)^ie_{n-i}q^i\Bigr)\;\cdot\;
    \Bigl(\sum_{i=0}^n(-1)^ie_iq^i\Bigr)\;,
  \end{gather*}
  whose coefficient at $q^1$ is \eqref{r1det} multiplied by $(-1)^n$.
  Now, substituting \eqref{r1det} in the numerator of
  \eqref{IW-form}, we get the identity
  \begin{equation}\label{preomid}
  \begin{gathered}
    \sum_{\substack{p\ge0\\i_1>\cdots>i_p\ge0}}\hspace{-2ex}
    (-1)^{p+i_1+\dots+i_p}
    s_{\((i_1+3),\dots,(i_p+3)\,|\,i_1,\dots,i_p\)}=
    \prod_{i\le j}^{n}(1-x_i x_j)\;\cdot\;
    \Bigl(\sum_je_j^2-e_{j+1}e_{j-1}\Bigr)=\\[-1.5ex]=
    \Bigl(\sum_{k=0}^{\frac{n(n+1)}2}(-1)^ke_k\comp h_2\Bigr)\;\cdot\;
    \Bigl(\sum_{j=0}^n e_j^2 - e_{j+1}e_{j-1}\Bigr)\;.
  \end{gathered}
  \end{equation}
  Since it holds for all $n$, we can consider \eqref{preomid} as an identity in the complete ring of symmetric functions (in the infinite set of variables) and apply to  \eqref{preomid} the $\omega$-involution, which exchanges $s_\l\leftrightarrow s_{\l'}$, $e_k\leftrightarrow h_k$. So, we get
  $$\sum_{\substack{p\ge0\\i_1>\cdots>i_p\ge0}}\hspace{-2ex}
     (-1)^{p+i_1+\dots+i_p}
     s_{\(i_1,\dots,i_p\,|\,(i_1+3),\dots,(i_p+3)\)}=
     \Bigl(\sum_k(-1)^ke_k\comp e_2 \Bigr)\;\cdot\;
     \Bigl(\sum_{j}h_j^2-h_{j+1}h_{j-1}\Bigr)\;,
  $$
  whose right hand side coincides with \eqref{KGrECpre}
\end{proof}

\begin{corollary}\label{doubledegree}
  For any decreasing sequence ${(N-3)\ge i_1>\dots>i_p>0}$ the
  irreducible representation
  $\pi_{\((i_1-1),\dots,(i_p-1)\,|\,( i_1+2),\dots,(i_p+2)\)}$
  appears in the both $\GL_N$-modules
  $$\L^{i_1+\dots+i_p}(\L^2W)\ox{A}_{(p)}\subset
    \OP_{k\ge 0}\L^{k}\(\L^2 W\)\ox A
  $$
  with multiplicity one and comes down into
  $R_{p,q}=H^p(K^\bdot_\Gr)_{(q)}$ with $q=p+i_1+\dots+i_p$\;.
\end{corollary}

\begin{proof}
  Using the Weyl formula\footnote{comp. with \cite[ch.\,1.5,\;example\,9(a)]{MD}} we can write $\chi_K$ as
  \begin{equation}\label{Wexpchi}
    \chi_{\L^\bdot(\L^2W[-1])}\cdot\chi_A=\hspace{-3ex}
    \sum_{\substack{p\ge 0\\N\ge i_1>\dots>i_p>0}}\hspace{-3ex}
    s_{\((i_1-1),\dots,(i_p-1)\,|\,i_1,\dots,i_p\)}\cdot
    \sum_{j\ge0}s_{[j,j]}\;.
  \end{equation}
  The classical Littlewood--Richardson rule\footnote{see \cite{Fu2}, \cite{MD}} implies that
  $s_{\((i_1-1),\dots,(i_p-1)\,|\,( i_1+2),\dots,(i_p+2)\)}$
  (staying in Lemma \ref{gr_char}) could appear in \eqref{Wexpchi}
  only as the product $s_{\((i_1-1),\dots,(i_p-1)\,|\,i_1,\dots,i_p\)}\cdot
  s_{[p,p]}$\,.
  Thus, the irreducible component
  $\pi_{\((i_1-1),\dots,(i_p-1)\,|\,( i_1+2),\dots,(i_p+2)\)}$
  comes into $H^p(K^\bdot_\Gr)_{(p+i_1+\dots+i_p)}$
  exactly from $\L^{i_1+\dots+i_p}(\L^2W)\ox{A}_{(p)}$, where it sits with multiplicity 1.
\end{proof}

\subsection{Dual Lie s-algebra}\label{G2Nstep2}
  Let us fix some coordinates $x^1,\dots,\,x^N$
  on $W$ and write $x^{ij}=x^i\wedge x^j$ for the
  corresponding basis of $L_1=\L^2W^*$, which is the generating
  space for the graded Lie s-algebra $L$ Koszul dual to $A$ in the
  sense of \ref{QDA}. Since
  $\Sym^2(\L^2W^*)\simeq\pi_{[2,2]}^*\oplus\pi_{[1^4]}^*$, it follows
  from \eqref{Ldef}, \eqref{Lbrac} that $L_2=\L^4W^*$ and the
  relations for $L$ have the form
  \begin{equation}\label{G2NLrel}
    [x^{ij},x^{k\ell}]=-[x^{ik},x^{j\ell}]=
    [x^{i\ell}, x^{jk}]\;.
  \end{equation}
  Let us write shortly $L'_n=L_n/\([\Lt,\Lt]\cap L_n\)=H^1(\Lt,\CC)_{(n)}=R_{n-1,n}$ for the space of degree $n$ generators in subalgebra $\Lt$.

\begin{lemma}\label{generators_lie_grasm}
  $L'_n\simeq\pi_{(n-2|n+1)}^*$ for each $n$ in range $2\le
  n\le N-2$.
\end{lemma}

\begin{proof}
  Induction on $n$. For $n=2$ we have $L'_{2}=L_2=\L^4W^*$
  generated by commutators \eqref{G2NLrel}. For $n>2$ we have
  the surjective commutator map
  $$L'_{n-1}\ox\L^2W^*
    \rOnto^{\quad a\ox x^{ij}\mapsto[a,x^{ij}]\quad}L'_n\;.
  $$
  By the inductive assumption and the Littlewood--Richardson rule,
  the left hand side has the following irreducible $\GL_N$-module decomposition
  \begin{equation}\label{commut_repr}
     \pi_{[1^2]}^*\ox\pi_{(n-3|n)}^*\simeq
     \pi_{(n-2,0|n,0)}^*\oplus
     \pi_{(n-2|n+1)}^* \oplus
     \pi_{(n-3,0|n,1)}^*\oplus
     \pi_{(n-3,0|n+1,0)}^*\oplus
     \pi_{(n-3|n+2)}^*\;.
   \end{equation}
   The s-Jacobi identity implies that
   $[a,[x^{ij},x^{k\ell}]]=[[a,x^{ij}],x^{k\ell}]+[[a,x^{k\ell}],x^{ij}]$\,.
   Since the left hand side here vanishes in $L'_n=\(\Lt/[\Lt,\Lt]\)_{(n)}$\,, the skew symmetrization operator taking
   $$\sum_\nu c_\nu\,[[a_\nu,x^{i_\nu j_\nu}],x^{k_\nu\ell_\nu}]
     \longmapsto\frac12\,\sum_\nu c_\nu\(
     [[a_\nu,x^{i_\nu j_\nu}],x^{k_\nu\ell_\nu}]-
     [[a_\nu,x^{k_\nu\ell_\nu}],x^{i_\nu j_\nu}]\)
   $$
   acts on $L'_n$ as the identity. On the other hand, it annihilates all the irreducible summands of \eqref{commut_repr} except for the second one, which definitely has to appear in $R_{n-1,n}=H^1(\Lt,\CC)_{(n)}=L'_n$ by the corollary \ref{doubledegree}.
\end{proof}

\subsection{Proof of Theorem \ref{syzygies}}\label{G2Nstep4}
  Since $\fA$ is (super) skew commutative and Koszul, its quadratic
  dual algebra $\fB=U(\fL)$ is an universal enveloping algebra for some
  graded Lie s-\-al\-geb\-ra $\fL$ such that $H^\bdot(\fL,\CC)=\fA$. It
  follows from Corollary \ref{doubledegree} and Lemma
  \ref{generators_lie_grasm} that there is a surjective homomorphism of
  associative algebras $H^{\bdot}(\Lt,\CC)\rOnto\fA$\,. Since 2-th Lie
  algebra cohomologies describe the relations, we conclude that the
  relations of $\Lt$ contain the ones of $\fL$. Because the both
  algebras are generated by the same vector space, there is a
  surjective $\GL$-module homomorphism $\fL\rOnto^\psi\Lt$. It follows
  from Corollary \ref{FLch} that $\chi_{\Lt}=\chi_{\fL}$, \ie the Euler
  $\GL$-cha\-rac\-ter of $\ker\psi$ vanishes. This forces $\ker\psi$ to
  have each irreducible  $\GL$-module equal number of times in its even
  and odd parts. But it is impossible, because the Young diagrams
  appearing in odd and even parts of $\fL$ have different number of
  cells modulo 4 (they consist of $4k+2$ and $4k$ cells respectively).
  Thus, $\ker\psi=0$ and theorem \ref{syzygies} is completely proven.

\section{Appendix: koszulity of hook algebras}

\noindent
  In this section we use the notations fixed in sec.\,\ref{GLnot}.

\subsection{Hook algebras}
  We call $\G$-shaped diagram $\G=(\a|\b)$ {\it a hook\/} of {\it
  width\/} $\a+1$ and {\it hight\/} $\b+1$. We say that two hooks
  $\G_1=(\a_1|\b_1)$, $\G_2=(\a_2|\b_2)$ are compatible and write
  $\G_1>\G_2$, if their union $\G_1\sqcup\G_2=(\a_1,\a_2|\b_1\b_2)$
  is a valid Young diagram, \ie has $\a_1>\a_2$ and $\b_1>\b_2$.

  Let us fix an ordered collection of pairwise compatible hooks
  \begin{equation}\label{hookcol}
     \G_1>\G_2>\;\cdots\;>\G_m
  \end{equation}
  equipped with some {\it internal parities\/} $|\G_i|\in\ZZ/(2)$.
  We assume that the heights of all hooks are bounded by $k$ and
  write $\pi_i=\pi_{\G_i}$ for the corresponding irreducible
  $\GL_k$-modules. For any increasing collection of indexes $I=\VEC
  i,s\subset[1\,.\,.\,m]$ we denote by
  \begin{equation}\label{hookYD}
     \G_I=\G_{\FAM i,s}=\G_{i_1}\sqcup\G_{i_2}\sqcup\cdots\sqcup\G_{i_s}
  \end{equation}
  the Young diagram build from the corresponding hooks and write
  $\pi_I=\pi_{\G_I}$ for the associated irreducible $\GL_k$-module.
  It follows from the Littlewood--Richardson rule that $\pi_I$
  appears in the irreducible decomposition of
  $\pi_{i_1}\ox\pi_{i_2}\ox\cdots\ox\pi_{i_s}$ with multiplicity
  one. Thus, there is a canonical $\GL_k$-equivariant projection
  \begin{equation}\label{hookmult}
     \mu_I:\pi_{i_1}\ox\pi_{i_2}\ox\cdots\ox\pi_{i_s}\rOnto\pi_I\;.
  \end{equation}
  Hence, $\GL_k$-module $\fA=\OP_I\pi_I$ admits a
  $\GL_k$-equivariant associative algebra structure whose
  multiplication satisfies the relations
  \begin{enumerate}
    \item
    $x\cdot y=(-1)^{|\G_i|\cdot|\G_j|}\,y\cdot x$ for any
    $x\in\pi_i$, $y\in\pi_j$ and any choice of $i\ne j$ in the
    range $[1\,.\,.\,m]$;
    \item
    $x\cdot y=0$ for any $x,y\in\pi_i$ and any choice of
    $i\in[1\,.\,.\,m]$;
    \item
    $x_1\cdot x_2\cdot\cdots\cdot x_s=\mu_I\VEC x,s$ for any choice
    of strictly increasing indexes
    $$I=\VEC i,s\subset[1\,.\,.\,m]
    $$
    and any collection of $x_\nu\in\pi_{i_\nu}$.
  \end{enumerate}
  In other words, $\fA$ is s-commutative \wrt the internal parity,
  it is generated by $\ZZ/(2)$-graded vector space
  $\fA_1=\OP_{i=1}^m\pi_{\G_i}$, and all non zero multiplication maps in $\fA$ are induced by projections \eqref{hookmult}. There is also  {\it an external grading\/}
  $$\fA=\OP_{s=0}^m\fA_s\;,\quad\text{where }
    \fA_s=\OP_{\#I=s}\pi_{\G_I}\;,\quad\fA_0=\CC\;.
  $$
  We call $\fA$ {\it a hook algebra\/} associated with hooks
  \eqref{hookcol} and write $\fA\VEC\G,m$, if the precise reference on
  the hooks is important.  In sec.\,\ref{hookproof}--\ref{PBWproof} we
  will prove

  \begin{theorem}\label{th:hookQK}
     Any hook algebra $\fA$ is quadratic and Koszul.
  \end{theorem}

  \noindent
  Since the algebra $\fA$ described in \eqref{G2NA1}--\eqref{eq:fAdef}
  is a hook algebra build from the hooks $\G_i=(i-1|i+2)$\,,
  $1\le i\le(N-3)$ of parities $|\G_i|=i\pmod 2$\,, we get

  \begin{corollary}\label{G2Nstep3}
     The algebra $\fA$, used in sec.\,\ref{answer_grasm},
     is quadratic and Koszul.
  \end{corollary}

\subsection{Proof of theorem \ref{th:hookQK}.}\label{hookproof}
  For each $\pi_i=\pi_{\G_i}$ we fix the standard basis
  labeled by the Young tableaux $T_{\G_i}$ of shape $\G_i$
  on the alphabet $[1\,.\,.\,k]$ (see sec.\,\ref{GLNchars}).

  Let $\wtd\fA$ be an algebra spanned by $\fA_1$ and satisfying only the first two sets of the relations for $\fA$, \ie $\wtd\fA$ is s-commutative \wrt the internal parity and satisfies $\pi_i\cdot\pi_i=0$ for all $i$. Thus, $\wtd\fA$ is a quadratic monomial algebra, in particular, it is automatically Koszul\footnote{see \cite{PP}}. With respect to the $\GL_k$-action, the graded components of $\wtd\fA$ are decomposed as
  $$\wtd\fA_s=\OP_{\#I=s}\pi_{i_1}\ox\pi_{i_2}\ox\cdots\ox\pi_{i_s}
  $$
  where the sum runs over all strictly increasing collections $I=\VEC i,s\subset[1\,.\,.\,m]$. Tensor products of the standard basic vectors from $\pi_i$ form a basis for $\wtd\fA_s$. We call these products {\it standard basic monomials\/}. They are numbered by Young diagrams \eqref{hookYD} filled by numbers from range $[1\,.\,.\,m]$ in such a way that each hook $\G_{i_\nu}\subset\G_I$ is a valid Young tableau but the whole $\G_I$ may be not. Let us call these filled diagrams {\it hooked tableaux\/} or \hqt x for shortness.

  We write $x_T\in\wtd\fA$ for the standard basic monomial corresponding to an \hqt\ $T$. Note that $x_S\cdot x_T=0$, if the underlying Young diagrams contain common hooks. Otherwise, $x_S\cdot x_T=\pm x_{S\cdot T}$, where $S\cdot T$ is build from $S$, $T$ by rearranging their hooks in strictly decreasing order. Thus, $x_T$ do actually behave as monomials.

  The hook algebra can be presented
  as $\fA=\wtd\fA/J$, where $J=\OP_sJ_s$ is a graded ideal whose
  components split \wrt the $\GL_k$-action as
  $$J_s=\OP_{\#I=s}J_I\;,\quad\text{where }
    J_I=
    \ker\(\pi_{i_1}\ox\pi_{i_2}\ox\cdots\ox\pi_{i_s}\rTo^{\mu_I}\pi_I\)\;.
  $$
  To show that $\fA$ is quadratic Koszul algebra, we will equip the
  set of all standard basic monomials $x_T\in\wtd\fA$ with a
  preorder $\preccurlyeq$ satisfying following two
  PBW-type\footnote{see \cite{PP} for non-commutative
  version of the Poincare--Birkhof--Witt theory} conditions:
  \begin{gather}
     \hbox{\parbox{.8\textwidth}
           {$x_S\preccurlyeq x_T\THEN x_{R\cdot S}\preccurlyeq x_{R\cdot T}$
            for any \hqt x $R$, $S$, $T$\,;
           }
          }\label{MCcond}
     \\[.5ex]
     \hbox{\parbox{.8\textwidth}
           {for any \hqt\ $T$ which is not a
            valid tableau there exist an element $h_T\in
            (J_2\cdot\wtd\fA)\cap J_I$ (uniquely determined by
            $T$) such that
            $$h_T-x_T=\hspace{-2ex}
              \sum\limits_{\substack{S\prec T\\I(S)=I(T)}}\hspace{-2ex}
              c_{S}x_{S}\quad\text{for some}\quad
              c_S\in\ZZ\;.
            $$
           }
          }\label{PBWcond}
  \end{gather}
  The second condition implies that each basic monomial $x_T$ is
  congruent modulo the elements $h_T$ to some monomial $x_{T'}$ whose
  \hqt\ $T'$ is a valid Young tableau. Since the images of the latter
  monomials form a basis for the vector space $\fA=\OP_I\pi_I=\wtd\fA/J$\,,
  we conclude that the elements $h_T$ generate $J$ as a vector
  space. Because $h_T$ lay in the ideal $(J_2)$ spanned by the
  quadratic component of $J$, we get $J=(J_2)$. Thus, $\fA$ is
  quadratic.

  Further, let $J^\circ$ be a monomial ideal spanned by the leading
  monomials of the elements from $J$. By the same reasons as above,
  $J^\circ$ is generated by all $x_T$ such that $T$ is not a valid
  tableau. Now the same arguments as in \cite[ch.\,3]{PP} show that
  koszulity of the monomial quadratic algebra
  $\fA^\circ=\wtd\fA/J^\circ$ implies the koszulity of $\fA$.

  Indeed, the multiplicative condition \eqref{MCcond} implies that
  the preorder in question induces a filtration on the bar complex
  of the hook algebra $\fA$ such that the associated graded complex
  is the bar complex of the Koszul algebra $\fA^\circ$. Computing
  $\ext_{\fA}(\CC,\CC)$ via the spectral sequence associated with
  this filtration, we get $\ext_{\fA^\circ}(\CC,\CC)$ as the first term
  of the sequence. It shows that $\ext_{\fA}^{i,j}(\CC,\CC)=0$ for
  $i\ne j$ (see details in \cite[ch.\,3]{PP}).

  Thus, to finish the proof of theorem \ref{th:hookQK}, it remains
  to equip the set of \hqt x with a preorder satisfying the
  PBW-properties \eqref{MCcond}--\eqref{PBWcond}.

\subsubsection{PBW-preorder on \hqt x}\label{hookord}
  Consider an \hqt\ $T$ whose Young diagram is the union of strictly
  decreasing hooks $\G_{i_1}>\G_{i_2}>\cdots>\G_{i_s}$. For any
  $\mu\in[1\,.\,.\,k]$, $\nu\in[1\,.\,.\,m]$ let $\chi^T(\mu,\nu)$ be a number of times the element $\mu$ does appear in the $\nu$-th hook
  $\G_{i_\nu}$ of $T$. For a fixed $\mu$ we consider the numbers
  $\chi^T(\mu,\nu)$ as the components of $m$-dimensional vector
  $$\xi_\mu^T=\(\chi^T(\mu,1),\chi^T(\mu,2),\ldots,\chi^T(\mu,m)\)
  $$
  where $\chi^T(\mu,\nu)=0$ when $\nu>s$. For example, if we deal
  with $\GL_4$-equivariant hook algebra built from $m=3$ hooks, then
  2-hooked tableau

  \centerline{$T=\hbox{\footnotesize$\young(1124,3122,43)$}$}

  \noindent
  produces four 3-component vectors
  \begin{equation}\label{chiexample}
  \begin{aligned}
     \chi_1^T&=(2,1,0)&\chi_3^T&=(1,1,0)\\
     \chi_2^T&=(1,2,0)&\chi_4^T&=(2,0,0)\;.
  \end{aligned}
  \end{equation}
  Note that a diagram $T$ is not uniquely defined by the
  collection of $k$ vectors
  $$\chi^T=\VEC\chi^T,k\;.
  $$
  For example, collection \eqref{chiexample} also comes from
  the \hqt

  \centerline{$T'=\hbox{\footnotesize$\young(1144,2123,32)$}$}

  \noindent
  and some others.

  We will compare the vectors $\chi^T_\mu$ using inverse right
  lexicographic ordering, \ie we say that
  $\(\chi^S_\mu(1),\chi^S_\mu(2),\dots,\chi^S_\mu(m)\)<
    \(\chi^T_\mu(1),\chi^T_\mu(2),\dots,\chi^T_\mu(m)\)$\,,
  if
  \begin{equation}\label{eq:irl-preorder}
    \chi^S(\mu,\nu_0)>\chi^T(\mu,\nu_0)\quad\&\quad
    \forall\;\nu>\nu_0\quad\chi^S(\mu,\nu)=\chi^T(\mu,\nu)\;.
  \end{equation}
  We say that $S\prec T$, if $\VEC\chi^S,k<\VEC\chi^T,k$ \wrt the
  inverse right lexicographic ordering, \ie if
  $\chi^S_\mu=\chi^T_\mu$ for all $\mu>\mu_0$ and
  $\chi^S_{\mu_0}>\chi^T_{\nu_0}$ in the sense \eqref{eq:irl-preorder}. By the definition, the condition $S\preccurlyeq T$
  means either the strong inequality $S\prec T$ or the coincidence
  $\VEC\chi^S,k=\VEC\chi^T,k$.

  Thus, the relation $\preccurlyeq$ gives a preorder on the set of \hqt x and two \hqt x are equivalent \wrt this preorder iff their fillings differ by a permutation preserving the content of each hook. This preorder evidently satisfies the multiplicative condition \eqref{MCcond}. It remains to construct special elements $h_T$ satisfying the  PBW-condition \eqref{PBWcond}.

\subsubsection{PBW-basis for $J$}\label{PBWproof}
     Consider an arbitrary hooked tableau $T$ which is not a valid
     tableau. Let the Young diagram of $T$ consist of
     hooks $\G_1>\G_2>\dots>\G_s$. We take minimal $i$ such that a
     hooked subtableau of $T$ formed by
     $\G_{i+1}\sqcup\G_{i+2}\sqcup\dots\sqcup\G_s$ is a valid
     tableau. Then a subtableau $D\subset T$ formed by
     $\G_i\sqcup\G_{i+1}$ is not valid. This can happen by two
     reasons (we write $\d_{p,q}$ for an element staying in $p$-th
     row and $q$-th column of $D$):
     \begin{enumerate}
        \item[(A)]
        $\exists\;k\ge2\;:\;\d_{k,2}<\d_{k,1}$, \ie the wrong
        inequality appears in some row of $D$ (but not in the first);
        \item[(B)]
        condition (A) fails but $\exists\;\ell,k\ge2\;:\;\d_{1,\ell}\ge\d_{k,\ell}$, \ie all rows of $D$ are valid and a wrong inequality appears in some column (but not in the first).
     \end{enumerate}
     We will construct $h_T$ separately for each case. However, in the both cases we construct $h_T$ as an element of the space
     \begin{gather}\label{eq:sJ2ker}
       \pi_1\ox\cdots\ox\pi_{i-1}\ox J_{i,i+1}\ox
       \pi_{i+2}\ox\cdots\ox\pi_s\;,\quad\text{where}\\\label{eq:J2ker}
       J_{i,i+1}=\ker\(\pi_{\G_i}\ox\pi_{\G_{i+1}}
       \rTo^{\quad\mu_{\G_i,\G_{i+1}}\quad}
       \pi_{\G_i\sqcup\G_{i+1}}\)\subset J_2\;.
     \end{gather}
     It follows from \LRr that the representation $\pi_T$, which is the target of the multiplication map \eqref{hookmult}, comes with multiplicity one in the space $\pi_1\ox\cdots\ox\pi_{i-1}\ox\pi_{\G_i\sqcup\G_{i+1}}\ox \pi_{i+2}\ox\cdots\ox\pi_s$\,, \ie the multiplication \eqref{hookmult} is factorized through this space and we can canonically include \eqref{eq:sJ2ker} into $J_s$. Further, we will present $h_T$ as $h_T=\pm x_S\cdot h_D\cdot x_R$, where $S$ and $R$ are formed by hooks $\G_\nu$ with $\nu<i$ and $\nu>i+1$ respectively, $D=\G_i\sqcup\G_{i+1}\subset T$ is the subdiagram formed by $i$-th and $(i+1)$-th consequent hooks, and $h_D$ lies in $J_{i,i+1}$ from \eqref{eq:J2ker}.

     Starting from this moment, we restrict ourself by this subdiagram $D$. Let its shape be $\d=(\a_1,\a_2|\b_1,\b_2)$, \ie $D=\G_1\sqcup\G_2$, where $\G_1=(\a_1|\b_1)$, $\G_2=(\a_2|\b_2)$.
     Consider the diagrams:
     $$\begin{aligned}
          \G_1'&=(0|\b_1)\quad
          \text{(\ie the first column of $\G_1$)}\;,\\
          \ovl\G_1&=(\a_1-1|0)=\G_1\smallsetminus\G_1'\;,\\
          \ovl\G_2&=(\a_1-1,\a_2-1|\b_2+1,0)=\ovl\G_1\sqcup\G_2
       \end{aligned}
     $$
     (\ie we split the first hook $\G_1=\G_1'\sqcup\ovl\G_1$ into
     the first column and remaining part of the row and form
     $\ovl\G_2$ by putting this row on top of $\G_2$). Recall that
     we write $\L^\l$ for the tensor product of exterior powers
     $\L^{\l_i'}V$ corresponding to the columns of a given Young
     diagram $\l$ (see \eqref{eq:La-space}). The factorization
     through the column exchange relations $\L^D\rOnto\pi_D$ can be
     considered as hauling down through the diagram
     \begin{equation}\label{eq:ER2dec}
     \begin{diagram}
        &&&&\L^D=\L^{\G_1'}\ox\L^{\ovl\G_2}\\
        &&&\ldTo<{1\ox\a}
            &&\rdTo>{1\ox\e}\\
        &&\L^{\G_1'}\ox\L^{\ovl\G_1}\ox\L^{\G_2}
         &&&&\L^{\G_1'}\ox\pi_{\ovl\G_2}\\
        &\ldEqto
         &&\rdTo>{1\ox\e\ox\e}
           &&\ldTo<{1\ox\r}\\
        \L^{\G_1}\ox\L^{\G_2}
        &&&&\pi_{\G_1'}\ox\pi_{\ovl\G_1}\ox\pi_{\G_2}
            &&\dEqto\\
        &\rdTo>{\e\ox\e}
         &&\ldTo<{\tau\ox1}
           &&\rdTo>{1\ox\s}\\
        &&\pi_{\G_1}\ox\pi_{\G_2}
          &&&&\pi_{\G_1'}\ox\pi_{\ovl\G_2}\\
        &&&\rdTo>\mu
           &&\ldTo<\eta\\
        &&&&\pi_D
       \end{diagram}
     \end{equation}
     where $\pi_{\G_1}\ox\pi_{\G_2}\rTo^\mu\pi_D$ is the
     multiplication map, $\e$ stays for particular factorizations
     through the column exchange relations, the map
     $\a:\L^{\ovl\G_2}\rTo\L^{\ovl\G_1}\ox\L^{\G_2}$ is the
     alternation in columns of $\ovl\G_2$, and the maps
     \begin{gather*}
     \r:\pi_{\ovl\G_2}\rInto\pi_{\ovl\G_1}\ox\pi_{\G_2}\;,\quad
     \s:\pi_{\ovl\G_1}\ox\pi_{\G_2}\rOnto\pi_{\ovl\G_2}\\
     \tau:\pi_{\G_1'}\ox\pi_{\ovl\G_1}\rOnto\pi_{\G_1}\;,\quad
     \eta:\pi_{\G_1'}\ox\pi_{\ovl\G_2}\rOnto\pi_D
     \end{gather*}
     are canonical projections onto and an inclusion of an irreducible
     submodule of multiplicity one.
     Indeed, it follows from \LRr that
     $\pi_{\ovl\G_2}$ has multiplicity one in the product
     $\pi_{\ovl\G_1}\ox\pi_{\G_2}$ and $\pi_{\G_1}$ has multiplicity
     one in the product $\pi_{\ovl\G_1'}\ox\pi_{\ovl\G_1}$. Similarly,
     $\pi_D$ has multiplicity one in the product
     $\L^{\G_1'}\ox\pi_{\ovl\G_1}\ox\pi_{\G_2}$ and this implies that
     the bottom rhombus of the diagram \eqref{eq:ER2dec} is commutative
     up to multiplication by non zero scalar factor. A straightforward
     computation (but quite improper for typesetting and too long for
     being reproduced here) shows that the composition
     $$\L^{\ovl\G_2}\rTo^\a\L^{\ovl\G_1}\ox\L^{\G_2}\rTo^{\e\ox\e}
       \pi_{\ovl\G_1}\ox\pi_{\G_2}
     $$
     annihilates all the column exchange relations in $\L^{\ovl\G_2}$, \ie the top rhombus in \eqref{eq:ER2dec} is commutative up to multiplication by non zero scalar factor as well. Thus the whole diagram \eqref{eq:ER2dec} is commutative up to rescales at the nodes. Now we are ready to describe the elements $h_D$.

     In case (B) we take the rightmost column with the wrong inequality $\d_{1,\ell}\ge\d_{k,\ell}$ and consider an element
     $$e^D+e^S\in
       \L^{\G_1}\ox\L^{\G_2}=
       \L^{\G_1'}\ox\L^{\ovl\G_1}\ox\L^{\G_2}\;,
     $$
     where $S\prec D$ is obtained from $D$ by transposing the entries $\d_{1,\ell}$, $\d_{k,\ell}$. Since the image of this element in
     $\pi_{\G_1'}\ox\pi_{\ovl\G_2}$ is zero, it follows from the commutativity of the bottom rhombus in \eqref{eq:ER2dec} that
     $$h_D=x_D+x_S=\e\ox\e(e^D+e^S)\in\pi_{\G_1}\ox\pi_{\G_2}
     $$
     lies in $\ker\mu\subset J_2$ as required.

     In case (A) we take the maximal $k$ such that $d_{k,2}<d_{k,1}$ (\ie the lowest row with the wrong order) and consider an exchange relation \eqref{eq:exchangerel} that exchanges first $k$ elements of the second column in $D$ with all $k$-element ordered subsets of the first column, \ie an element
     $$\wht h_D=e^D-\sum_\s e^{\s D}=e^D-\sum_S a_Se^S\in
       \L^D=\L^{\G_1'}\ox\L^{\ovl\G_2}\;,
     $$
     where $S\prec D$ are obtained from $D$ by all the exchanges $\s$ in question. Since the image of $\wht h_D$ in $\pi_D$ is zero, the class of an element $1\ox\a(\wht h_D)\in\L^{\G_1}\ox\L^{\G_2}$ in the factor $\pi_{\G_1}\ox\pi_{\G_2}$ belongs to $\ker\mu\subset J_2$. At the same time the difference
     $$\Bigl(x_D-\sum_S a_Sx_S\Bigr)-\e\ox\e\bigl(1\ox\a(\wht h_D)\bigr)\in
       \pi_{\G_1}\ox\pi_{\G_2}
     $$
     is a sum of $\e\ox\e$-images of elements from $\L^{\G_1}\ox\L^{\G_2}$ that have a form considered above in the case (B). In particular, this difference lies in $J_2$ as well. We conclude that
     $$h_T=x_D-\sum_S a_Sx_S
     $$
     satisfies the required properties.


\pagebreak[3]

\end{document}